  %
  %



  \font \bbfive = bbm5
  \font \bbseven = bbm7
  \font \bbten = bbm10
  \font \eightbf = cmbx8
  \font \eighti = cmmi8 \skewchar \eighti = '177
  \font \eightit = cmti8
  \font \eightrm = cmr8
  \font \eightsl = cmsl8
  \font \eightsy = cmsy8 \skewchar \eightsy = '60
  \font \eighttt = cmtt8 \hyphenchar\eighttt = -1
  \font \msbm = msbm10
  \font \sixbf = cmbx6
  \font \sixi = cmmi6 \skewchar \sixi = '177
  \font \sixrm = cmr6
  \font \sixsy = cmsy6 \skewchar \sixsy = '60
  \font \tensc = cmcsc10
  
  \font \titlefont = cmr7 scaled \magstep4
  \scriptfont \bffam = \bbseven
  \scriptscriptfont \bffam = \bbfive
  \textfont \bffam = \bbten


  \def \Headlines #1#2{\nopagenumbers
    \voffset = 2\baselineskip
    \advance \vsize by -\voffset
    \headline {\ifnum \pageno = 1 \hfil
    \else \ifodd \pageno \tensc \hfil \lowercase {#1} \hfil \folio
    \else \tensc \folio \hfil \lowercase {#2} \hfil
    \fi \fi }}

  \def \Title #1{\vbox{\baselineskip 20pt \titlefont \noindent #1}}

  \def \Date #1 {\footnote {}{\eightit Date: #1.}}

  \def \Authors #1{\bigskip \bigskip \noindent #1}

  \long \def \Addresses #1{\begingroup \eightpoint \parindent0pt
\medskip #1\par \par \endgroup }

  \long \def \Abstract #1{\begingroup \eightpoint
  \bigskip \bigskip \noindent
  {\sc ABSTRACT.} #1\par \par \endgroup }

  \newcount \ccnt \ccnt = 0
  \newdimen \ldim
  \newdimen \rdim
  \newdimen \ndim
  \newdimen \mdim
  \ldim = 0.1\hsize
  \rdim = 0.1\hsize
  \ndim = 15pt
  \mdim = \hsize
  \advance \mdim by -\ldim
  \advance \mdim by -\rdim
  \advance \mdim by -\ndim
  \def \triplebox #1#2#3{\noindent \hbox {\hbox to \ldim {\hfil
#1}\hbox to \mdim { #2\dotfill }\hbox to \rdim {\hbox to \ndim {\hfil
#3}\hfil }}\par }

  \def \pag #1 ... #2{\triplebox {\global \advance \ccnt by 1 \number
\ccnt .}{#2}{#1}}
  \def \ppg #1 ... #2{\triplebox {}{#2}{#1}}


  \def \vg #1{\ifx #1\null \null \else
    \ifx #1\ { }\else
    \ifx #1,,\else
    \ifx #1..\else
    \ifx #1;;\else
    \ifx #1::\else
    \ifx #1''\else
    \ifx #1--\else
    \ifx #1))\else
    { }#1\fi \fi \fi \fi \fi \fi \fi \fi \fi }

  \def \goodbreak {\vskip0pt plus.1\vsize \penalty -250 \vskip0pt
plus-.1\vsize }

  \newcount \secno \secno = 0
  \newcount \stno

  \def \seqnumbering {\global \advance \stno by 1
    \number \secno .\number \stno }

  \def \label #1{\def\localsystemvariable {\number \secno
    \ifnum \number \stno = 0\else .\number \stno \fi }\global \edef
    #1{\localsystemvariable }}
  \def\ilabel #1{\global \edef #1{\numbertoi}\ignorespaces}

  \def\section #1{\global\def\SectionName{#1}\stno = 0 \global
\advance \secno by 1 \bigskip \bigskip \goodbreak \noindent {\bf
\number \secno .\enspace #1.}\medskip \noindent \ignorespaces}

  \long \def \sysstate #1#2#3{\medbreak \noindent {\bf
\seqnumbering .\enspace #1.\enspace }{#2#3\vskip 0pt}\medbreak }
  \def \state #1 #2\par {\sysstate {#1}{\sl }{#2}}
  \def \definition #1\par {\sysstate {Definition}{\rm }{#1}}


  \def \proof {\medbreak \noindent {\it Proof.\enspace }}
  \def \proofend {\ifmmode \eqno \square \else \hfill \square
\looseness = -1 \medbreak \fi }

  \def \$#1{#1 $$$$ #1}

  \newcount \zitemno \zitemno = 0
  \def\numbertoi {{\rm(\ifcase \zitemno \or i\or ii\or iii\or iv\or
v\or vi\or vii\or viii\or ix\or x\or xi\or xii\fi)}}
  \def\zitemcntr {\global \advance \zitemno by 1 \numbertoi}
  \def\iItem #1{\medskip \item {#1}}
  \def\Item #1{\smallskip \item {#1}}
  \def\izitem {\zitemno = 0\iItem {\zitemcntr}}
  \def\zitem {\Item {\zitemcntr}}

  \newcount \footno \footno = 1
  \newcount \halffootno \footno = 1
  \def\footcntr {\global \advance \footno by 1
  \halffootno =\footno
  \divide \halffootno by 2
  $^{\number\halffootno}$}


  \def \({\left (}
  \def \){\right )}
  \def \[{\left \Vert }
  \def \]{\right \Vert }
  \def \*{\otimes }
  \def \+{\oplus }
  \def \:{\colon }
  \def \<{\left \langle }
  \def \>{\right \rangle }
  \def \text #1{\hbox {\rm #1}}
  \def \and {\hbox {,\quad and \quad }}
  
  \def \calcat #1{\,{\vrule height8pt depth4pt}_{\,#1}}

  \def \crossproduct {{\hbox {\msbm o}}}
  \def \cstar {$C^*$}
  \def \for #1{,\quad #1}
  \def \inv {^{-1}}
  
  \def \square {\hbox {$\sqcap \!\!\!\!\sqcup $}}
  \def \stress #1{{\it #1}\/}
  
  \def \x {\times }
  \def \|{\Vert }
  \def \inv {^{-1}}

  \newskip \ttglue

  \def \eightpoint {\def \rm {\fam0 \eightrm }%
  \textfont0 = \eightrm
  \scriptfont0 = \sixrm \scriptscriptfont0 = \fiverm
  \textfont1 = \eighti
  \scriptfont1 = \sixi \scriptscriptfont1 = \fivei
  \textfont2 = \eightsy
  \scriptfont2 = \sixsy \scriptscriptfont2 = \fivesy
  \textfont3 = \tenex
  \scriptfont3 = \tenex \scriptscriptfont3 = \tenex
  \def \it {\fam \itfam \eightit }%
  \textfont \itfam = \eightit
  \def \sl {\fam \slfam \eightsl }%
  \textfont \slfam = \eightsl
  \def \bf {\fam \bffam \eightbf }%
  \textfont \bffam = \eightbf
  \scriptfont \bffam = \sixbf
  \scriptscriptfont \bffam = \fivebf
  \def \tt {\fam \ttfam \eighttt }%
  \textfont \ttfam = \eighttt
  \tt \ttglue = .5em plus.25em minus.15em
  \normalbaselineskip = 9pt
  \def \MF {{\manual opqr}\-{\manual stuq}}%
  \let \sc = \sixrm
  \let \big = \eightbig
  \setbox \strutbox = \hbox {\vrule height7pt depth2pt width0pt}%
  \normalbaselines \rm }


  \def\cite #1{{\rm [\bf #1\rm ]}}
  \def\scite #1#2{\cite{#1{\rm \hskip 0.7pt:\hskip 2pt #2}}}
  \def\lcite #1{#1}
  \def\fcite #1#2{#1}
  \def\bibitem#1#2#3#4{\smallskip \item {[#1]} #2, ``#3'', #4.}

  \def \references {
    \begingroup
    \bigskip \bigskip \goodbreak
    \eightpoint
    \centerline {\tensc References}
    \nobreak \medskip \frenchspacing }


  \def\uint{\int^{u}}
  \def\sint{\int^{su}}
  \def\C{{\bf C}}
  \def\Z{{\bf Z}}
  \def\B{{\cal B}}
  \def\Ft{{\cal F}}
  \def\H{{\cal H}}
  \def\J{{\cal J}}
  \def\K{{\cal K}}
  \def\W{{\cal W}}
  \def\M{{\cal M}}
  \def\N{{\cal N}}
  \def\P{{\cal P}}
  \def\Ma{{\cal M}_\a}
  \def\Na{{\cal N}_\a}
  \def\Pa{{\cal P}_\a}
  \def\I{{\cal I}}
  \def\O{{\cal O}}
  \def\J{{\cal J}}
  \def\Xa{{\cal X}_\a}
  \def\X{{\cal X}}
  \def\T{X}
  \def\Mult{M}
  \def\a{\alpha }
  \def\Gdual{{\^G}}
  \def\ah{\widehat{\alpha }}
  \def\d{\,d}
  \def\e#1{e_{#1}}
  \def\^#1{\widehat{#1}}
  \def\rcp{A\crossproduct_{\a,r} G}
  \def\fullcp{A\crossproduct_\a G}
  \def\V{V}
  \def\conv{\hbox{$*$}}
  \def\du#1#2{\,\hbox{$<$}#1,#2\hbox{$>$}\,}
  \def\tdu#1#2{\du{#1}{#2}}
  \def\idu#1#2{\,\overline{\hbox{$<$}#1,#2\hbox{$>$}}\,}
  \def\F#1#2{E_{#2}(#1)}
  \def\E#1{\F{#1}e}
  \def\Rip#1#2{\<#1,#2\>_R}
  \def\Lip#1#2{\<#1,#2\>_L}
  \def\char#1{\chi_{\null_{#1}}\!}
  \def\=#1{\buildrel (#1) \over =}
  \def\span{\text{span}}
  \def\rc{\buildrel{rc} \over \sim}
  \def\Alg#1{{Alg}(#1,\a)}
  \def\bigarrow{{\hbox to 1cm{\rightarrowfill}}}
  \def\A{C_0(\T)}
  \def\Dxyz{\Delta_{x,y}^{z}}
  \def\MR{Morita--Rieffel}
  \def\GFPA{generalized fixed point algebra}


  \def\Newpim{E1}
  \def\unconditional{E2}
  \def\TPA{E3}
  \def\Amena{E4}
  \def\Rep{E5}
  \def\FD{FD}
  \def\Green{G}
  \def\Halmos{H}
  \def\HR{HR}
  \def\Jensen{JT}
  \def\KR{KR}
  \def\McClanahan{M}
  \def\Ng{N}
  \def\Ped{P}
  \def\RieffelInduced{R1}
  \def\RieffelOA{R2}
  \def\RieffelTransf{R3}
  \def\RieffelProper{R4}
  \def\RieffelIntegrable{R5}
  \def\Yosida{Y}
  \def\Zettl{Z}


  \Headlines
  {Morita--Rieffel Equivalence and Spectral Theory} 
  {Ruy Exel} 

  \Date {April 19, 1999}

  \Title
  {Morita--Rieffel Equivalence and Spectral Theory for\hfill\break
  Integrable Automorphism Groups of C*-Algebras}

  \Authors
  {Ruy Exel\footnote {*}{\eightrm Partially supported by CNPq.}
  {\eightrm(exel@mtm.ufsc.br)}}

  \Addresses
  {Departamento de Matem\'atica;
  Universidade Federal de Santa Catarina;
  88010-970 Florian\'opolis SC;
  BRASIL.}

  \Abstract
  {Given a \cstar-dynamical system $(A,G,\a)$, we discuss conditions
under which subalgebras of the multiplier algebra $\Mult(A)$
consisting of fixed points for $\a$ are {\MR} equivalent to ideals in
the crossed product of $A$ by $G$.  In case $G$ is abelian we also
develop a spectral theory, giving a necessary and sufficient condition
for $\a$ to be equivalent to the dual action on the cross-sectional
\cstar-algebra of a Fell bundle.
  In our main application we show that a proper action of an abelian
group on a locally compact space is equivalent to a dual action.
  }

  \newcount \ccnt \ccnt = 0
  \def\.{\hbox to 70pt{\hfill}}
  \def\indx#1{\noindent\. #1 \dotfill \global \advance \ccnt by 1 \
\number \ccnt\.\par}
  \bigskip
  \centerline {CONTENTS}
  \medskip
  \indx {Introduction}
  \indx {Preliminaries}
  \indx {Integrable elements}
  \indx {The left structure of $\Xa$}
  \indx {Integral and Laurent operators}
  \indx {Abelian groups and the Fourier transform}
  \indx {The dual unitary group}
  \indx {Relative continuity}
  \indx {{\MR} equivalence}
  \indx {Dual action on Fell bundles}
  \indx {Spectral theory}
  \indx {Classical dynamical systems}
  \indx {An example}

   \section {Introduction}
  The object under study in the present work is a \cstar-dynamical
system $(A,G,\a)$, that is, a strongly continuous action $\a$ of a
locally compact topological group $G$ on a \cstar-algebra $A$.  The
questions we discuss fall broadly in two categories, the first one
being Rieffel's project, initiated in \cite{\RieffelProper} and
recently continued in \cite{\RieffelIntegrable}, of defining a
{\GFPA}, under suitable hypothesis, and proving it to be
  {\MR}\footnote{\footcntr}
  {\eightpoint Given the enormous contribution made by Rieffel on the
subject of Morita equivalence in the context of \cstar-algebras, I
believe that the term ``strong Morita equivalence'' should be
re-coined ``{\MR} equivalence''.  In fact I was told that similar
terms were already employed in the 1998 Great Plains Operator Theory
Symposium (GPOTS) held at Kansas State University.}
  equivalent to an ideal in the reduced crossed product algebra.  The
second issue we discuss is \stress{spectral theory} for $G$ abelian.

To describe what we mean by spectral theory suppose that $G$ is not
only abelian but also compact.  In this very simple case it is easy to
show (see e.g. \scite{\Newpim}{2.5}) that $A$ may be decomposed as the
closure of the direct sum $\oplus_{x\in\Gdual} B_x$, where for each
$x$ in the Pontryagin dual $\Gdual$ of $G$, the \stress{spectral
subspace} $B_x$ is the closed subspace of $A$ formed by the elements
$a\in A$ such that $\a_t(a)=\idu xt a$, for all $t$ in $G$.  As usual
we denote by $\du xt$ the duality between $\Gdual$ and $G$.

By spectral theory we mean a theory encompassing generalizations of
the above decomposition including situations where $G$ is not compact.
Returning for a moment to the compact case, one may prove that the
spectral subspaces $B_x$ defined above satisfy
  $$
  B_x B_y \subseteq B_{xy}
  \and
  B_x^* = B_{x\inv},
  $$
  for all $x$ and $y$ in $\Gdual$.  It therefore follows that the
collection $\B=\{B_x\}_{x\in\Gdual}$ forms a Fell bundle over the
discrete group $\Gdual$ (see \cite{\FD} for a comprehensive treatment
of the theory of Fell bundles, also referred to as \cstar-algebraic
bundles).  One may then prove \scite{\Amena}{4.7} that $A$ is
isomorphic to $C^*(\B)$, the cross-sectional \cstar-algebra of $\B$
\scite{\FD}{VIII.17.2}.  Moreover the natural isomorphism between $A$
and $C^*(\B)$ turns out to be covariant for the \stress{dual action}
\scite{\unconditional}{Section 5} of $G$ on $C^*(\B)$.  In short,
every compact abelian action is equivalent to a dual action.

If $G$ is abelian but no longer compact one can still consider the
dual action of $G$ on the cross-sectional \cstar-algebra of any given
Fell bundle over the (no longer discrete) dual group $\Gdual$.  We
therefore take this to be the model action for spectral theory and
hence when a given action is shown to be equivalent to a dual action
we shall consider that the goal of spectral theory has been achieved.

One of the reasons one might benefit from proving an action to be a
dual action is that Fell bundles can be classified
\scite{\TPA}{Theorem 7.3} up to stable equivalence by a twisted
partial action of the base group on the unit fiber algebra.  The
successful completion of this program is therefore likely to provide a
deep insight on the behavior of the given dynamical system.

In one of our main results, namely Corollary
\fcite{11.15}{CharacterizationOfDualaction}, we give a necessary and
sufficient condition for an action to be equivalent to a dual action,
thus solving a problem considered in \cite{\unconditional}.  According
to Theorem 5.5 in \cite{\unconditional}, in the case of a dual action
$A$ contains a dense set of $\a$-\stress{integrable elements}: an
element $a\in A$ is said to be \stress{$\a$-integrable}
\cite{\unconditional} (see also \cite{\RieffelIntegrable}) if, for
every $b\in A$, the functions
  $
  t \mapsto \a_t(a)b
  $
  and,
  $
  t \mapsto b\a_t(a)
  $
  are unconditionally integrable in the sense that the Bochner
integrals over compact subsets of $G$ form a converging net
(Definition \fcite{2.1}{AltDefinitionOfIntegr} below).  Accordingly,
the necessary and sufficient conditions of
\fcite{11.15}{CharacterizationOfDualaction} are expressed in terms of
the existence of certain sets of $\a$-integrable elements.

We are then able to verify that these conditions hold for a proper
action on a locally compact topological space and hence we obtain
(Corollary \fcite{12.5}{SpectralTheoremCommutativeCase}) that proper
actions are equivalent to dual actions.

At first glance it may seem that the questions related to the {\GFPA}
and {\MR} equivalence bear no relation to spectral theory but these
turn out to be intimately related subjects.  Among the indications
that this is so is Theorem \fcite{10.6}{MoritaForFellBundles} below,
according to which Rieffel's project may be satisfactorily carried out
when the action in question is a dual action.  In particular, the role
of the \stress{{\GFPA}} is played by the unit fiber algebra which
appears as a subalgebra of the multiplier algebra and is fixed by the
dual action.

At the root of the relationship between Rieffel's project mentioned
above and spectral theory is the notion of \stress{relative
continuity} which we would now like to describe in some detail.
Following \cite{\unconditional} and \cite{\RieffelIntegrable}
(although we sometimes use different notation and terminology) we let
$\Pa$ be the set of positive $\a$-integrable elements and $\Na =
\{a\in A : a^*a \in \Pa\}$.  It turns out that $\Pa$ is a hereditary
cone \scite{\unconditional}{6.6} which implies that $\Na$ is a left
ideal in $A$.  Moreover, it can be shown without much difficulty that
$\Na$ is a right Hilbert $\Mult_e(A)$--module, where $\Mult_e(A)$
refers to the subalgebra of the multiplier algebra $\Mult(A)$
consisting of the fixed points for $\a$.  The $\Mult_e(A)$--valued
inner-product $\Rip ab$ is given for any $a,b\in\Na$ by the
strict-unconditional integral of $\a_t(a^*b)$, with respect to $t$.

It is implicit in Rieffel's early work on this subject
\cite{\RieffelProper} that $\Na$, with its Hilbert module structure,
hides the clue to the {\MR} equivalence between the {\GFPA}, which has
yet to be defined, and ideals in the crossed product algebra.

Our starting point is a certain representation of the Hilbert module
$\Na$ as bounded operators from $\H$ to $L_2(G,\H)$, where $\H$ is any
fixed representation space for $A$ (and hence also for $\Mult(A)$).
Precisely speaking, we define a linear map
  $$
  \zeta : \Na \to \B(\H,L_2(G,\H))
  $$
  which satisfies
  $\zeta(a)m= \zeta(am)$, and
  $\Rip{a}{b} = \zeta(a)^*\zeta(b)$
  for all $a,b\in\Na$ and $m\in\Mult_e(A)$.
  Once $\Na$ is concretely represented, our ability to answer
questions related to it greatly increase.  In particular it is
immediate that the algebra of \stress{generalized compact operators}
  \cite{\RieffelInduced}
  is the closed linear span of the set
  $\{\zeta(a)\zeta(b)^* : a,b\in \Na\}$, which is a subset of
$\B(L_2(G,\H))$.

Since the reduced crossed product $\rcp$ is also an algebra of
operators on $L_2(G,\H)$, it makes sense to ask whether or not
  $$
  \zeta(a)\zeta(b)^*\in \rcp
  \leqno{(\seqnumbering)}
  \label \MemberShip
  $$
  for a given pair of elements $a,b\in\Na$.  We show in
\fcite{13.5}{RCForRankOneOperators} that in an important example this
fails more often than it holds.  When $G$ is abelian we give a
necessary and sufficient condition for \lcite{\MemberShip} to hold in
terms of the \stress{Fourier coefficients} $\F px$ of the element
$p:=a^*a$, defined by
  $$
  \F px = \sint_G \du xt \a_t(p) \d t
  \for x\in \Gdual,
  $$
  and that of $q:=b^*b$.
  Precisely, we show in Theorem \fcite{7.5}{IncredibleFormula} that
\lcite{\MemberShip} holds if and only if
  $$
  \lim_{z\to e} \[\F p{xz} \F qy - \F px \F q{zy}\] =0,
  \leqno{(\seqnumbering)}
  \label \IFONE
  $$
  uniformly in $x$ and $y$ in $\Gdual$.  A first consequence is that
\lcite{\MemberShip} depends only on the absolute value of $a$ and $b$.
Secondly it brings to light a relation between $\a$-integrable
elements: given any pair of $\a$-integrable elements $p$ and $q$ we
say that $p$ and $q$ are \stress{relatively continuous}, and denote it
as $p\rc q$, if \lcite{\IFONE} holds.

The relation ``\ $\rc$\ '' is not reflexive, transitive, or symmetric
(see \fcite{13.5}{RCForRankOneOperators}) but it enjoys some
surprising properties (see \fcite{8.3}{RCProperties},
\fcite{8.4}{Herege}, and \fcite{8.5}{IntIsRelContinuous}) such as
  $$
  0\leq a\leq b \rc c \Rightarrow a\rc c,
  $$
  whenever $a$, $b$, and $c$ are $\a$-integrable.

If one is to find a submodule $\N$ of $\Na$ which serves as the
imprimitivity bimodule for a {\MR} equivalence involving a subalgebra
of $\rcp$, then it would be convenient if $\N$ satisfied
$\zeta(\N)\zeta(\N)^*\subseteq\rcp$.  Therefore, for all $a,b\in\N$ it
must be that
  $a^*a\rc b^*b$ and hence, by the polarization formula, that
  $p\rc q$ for any $p,q\in\N^*\N$.

We then see that $\N^*\N$ must be a \stress{relatively continuous set}
in the sense that its elements are mutually relatively continuous.
The role of relative continuity in relation to {\MR} equivalence thus
becomes clear.  In particular one of the main hypothesis in Theorem
\fcite{9.2}{AbstractMR}, our main result related to {\MR} equivalence,
is that a certain set is relatively continuous.
  With respect to Theorem \fcite{11.14}{SpectralTheorem}, one of our main
results in spectral theory, relatively continuous sets also play a
crucial role.  In particular the continuity of the norm and
multiplication on the Fell bundle we construct there are derived from
the relatively continuous set $\W$ in the hypothesis of
\fcite{11.14}{SpectralTheorem}.

Relatively continuous sets are therefore crucial for our arguments to
be carried out in both fronts.  However, it could be argued that
requiring relative continuity beforehand diminishes the scope of
application of our main results.
  It would therefore be highly desirable to produce large relatively
continuous sets as a consequence of the existence of large (say dense)
sets of $\a$-integrable elements (see questions
  \fcite{9.4}{DensityOfRCCone},
  \fcite{9.5}{QuestionOnUFA},
  and
  \fcite{11.16}{MainQuestion}).
  However, except for Section \fcite{12}{ClassicalDynSystemsSection},
which deals with commutative algebras, we have nothing to offer in
this respect.

As the main application of our results we discuss in Section
\fcite{12}{ClassicalDynSystemsSection} the case of a classical
dynamical system
  $(\A,G,\a)$ based on an abelian group, therefore assuming as much
commutativity as possible.

By a result of Rieffel \scite{\RieffelIntegrable}{4.7} the linear span
of the positive $\a$-integrable elements is dense in $\A$ if and only
if $\a$ is a proper action in the usual sense.  Therefore it seems
reasonable to restrict our attention to proper actions.  In this case
it turns out that every $f$ in $C_c(\T)$ is relatively continuous with
respect to any other $\a$-integrable element (Proposition 12.3) and
hence it deserves to be called \stress{absolutely continuous}.  This
implies that $C_c(\T)$ is a relatively continuous set and hence we are
able to apply the general theory developed in the previous sections to
show the already mentioned fact  that proper actions are equivalent to
dual actions.

Given the essential role played by absolutely continuous elements in
our treatment of classical dynamical systems one could suspect them to
play a similar role in other situations.  However, in the case of the
action of\/ $\Z$ on the algebra $\K$ of compact operators given by
conjugation by the powers of the bilateral shift, we prove in
\fcite{13.6}{NoPositiveAbsolutelyContinuous} that the only positive
absolutely continuous element is the zero operator.  So it is not
clear, in general, how to produce a canonical relatively continuous
set.

This explains, in retrospect, why is it so important to postulate the
existence of the dense subalgebra $A_0$ in Definition 1.2 of
\cite{\RieffelProper}.  Since relative continuity is a
\stress{relative} relation, there seems not to be a universal
construction of $A_0$ that works.  Exploring the example of the action
of $\Z$ on $\K$, mentioned above, we are able to give a precise reason
why this is so by showing that a maximal relatively continuous cone is
not unique.

It is noteworthy that several tools from the Harmonic Analysis of
abelian groups are used in this work.  While this is quite natural
with respect to spectral theory, it is not so much so in relation to
the aspects of {\MR} equivalence that we discuss.  In fact, the very
statement of our main spectral theory result (Theorem 11.14) makes no
sense if $G$ is not commutative, given the use of the Pontryagin dual.
However, Rieffel's project of generalizing the concept of fixed point
algebra is meaningful for non-commutative groups and in fact many such
examples have already been treated \cite{\RieffelProper}.  This
indicates either that our option to use tools from Harmonic Analysis
is a poor one, or that the stakes are higher then previously thought!

  \section {Preliminaries}
  \label \BasicFactsSection
  Throughout this article we will fix a \cstar-algebra $A$ and a
strongly continuous action
  $$
  \a : G \to {\rm Aut}(A)
  $$
  of a locally compact group $G$ on $A$.  We will also suppose,
without loss of generality, that $A$ is a non-degenerate
\cstar-algebra of operators on a Hilbert space $\H$.

  Recall from \scite{\unconditional}{2.3} that a function
  $f\colon G\to A$ is said to be \stress{unconditionally integrable}
if it is Bochner integrable \scite{\Yosida}{V.5} over every relatively
compact subset $K\subseteq G$, and the net
  $$
  \left\{\int_K f(t)\d t\right\}_{K\in \K}
  $$
  converges in the norm topology of $A$, where $\K$ is the directed
set of all relatively compact subsets of $G$ ordered by inclusion.

Although this is not relevant for us here we should remark that it is
possible to prove that the notion of unconditional integrability
coincides with that of Pettis integrability at least for bounded
continuous functions.

If $f$ is unconditionally integrable we let
  $$
  \uint_G f(t) \d t :=
  \lim_K \int_K f(t)\d t.
  $$
  The superscript ``$u$'' in the integral above is meant to remind us
that we are speaking of an \stress{unconditional} integral.
  According to \scite{\unconditional}{6.1} we say that an element $a$
in $A$ is \stress{$\a$-integrable} if for every $b$ in $A$ one has
that the functions
  $$
  t\in G \mapsto \a_t(a)b \in A
  \and
  t\in G \mapsto b\a_t(a) \in A
  $$
  are unconditionally integrable (see also
\cite{\RieffelIntegrable}). In this case we denote by
  $$
   \sint_G \a_t(a) \d t
  $$
  the multiplier $(L,R)$ of $A$ given by
  $$
  L(b) = \uint_G \a_t(a) b \d t
  \and
  R(b) = \uint_G b \a_t(a) \d t,
  $$
  for every
  $b\in A$.  The superscript ``$su$'' standing for \stress{strict
unconditional} integral.

  For simplicity, in order to deal with left and right multiplication
at the same time, we will often adopt the following alternative
definition of $\a$-integrability:

  \definition
  \label \AltDefinitionOfIntegr
  An element $a\in A$ is said to be $\a$-integrable if, given any pair
$b,c$ of elements of the multiplier algebra $\Mult(A)$ of $A$ such
that either $b\in A$ and $c=1$, or $b=1$ and $c\in A$, we have that
the map
  $
  t\in G \mapsto b\a_t(a)c \in A
  $
  is unconditionally integrable (compare
\scite{\RieffelIntegrable}{2.2}).

  \sysstate{Remark}{\rm}
  {\label \RemarkOnDifferentTerminology We should remark that our
notion of \stress{$\a$-integrable} elements, taken from
\scite{\unconditional}{6.1}, is related to Rieffel's notion of
\stress{$\a$-proper} elements \scite{\RieffelIntegrable}{4.1} as
opposed to what is called \stress{order-integrable} in
\scite{\RieffelIntegrable}{1.4}.  See also
\scite{\RieffelIntegrable}{4.4}.}

  Throughout this work we shall adopt the following:

  \sysstate{Notation}{\rm}
  {\label\NotationMNP
  \izitem $\Pa = \{a\in A_+ : \text{ $a$ is $\a$-integrable}\}$,
  \zitem $\Na = \{a\in A : a^*a\in \Pa\}$,
  \zitem $\Ma = \span(\Pa)$.}

  We will denote also by $\a$ the usual (not necessarily strongly
continuous) extension of the action $\a$ to $\Mult(A)$.  The subset of
$\Mult(A)$ formed by the fixed points under $\a$ will be denoted
$\Mult_e(A)$.  It is easy to see that for any $\a$-integrable element
$a$ one has that $\sint_G\a_t(a)\d t\in \Mult_e(A)$.

  \state Proposition
  \label \BasicPropertiesOfModule
  \izitem $\Pa$ is a hereditary cone,
  \zitem \ilabel{\zIntegr} $\Ma$ consists of $\a$-integrable elements,
  \zitem \ilabel{\zaNinN} $\Pa$, $\Na$, and $\Ma$ are invariant under
$\a$,
  \zitem \ilabel{\zNMinN} $\Na\Mult_e(A)\subseteq\Na$,

  \proof
  That $\Pa$ is hereditary is precisely the content of
  \scite{\unconditional}{6.6}.  See also
\scite{\RieffelIntegrable}{2.7}.  It is obvious that a linear
combination of $\a$-integrable elements is again $\a$-integrable and
hence (ii) holds.

  Let $a$ be any $\a$-integrable element and let $s\in G$. Then for
every $K\in\K$ and $b,c$ as in \lcite{\AltDefinitionOfIntegr} we have
that
  $$
  \int_K b\a_t(\a_s(a))c \d t =
  \int_K b\a_{ts}(a)c \d t =
  \Delta(s\inv) \int_{Ks} b\a_t(a)c \d t,
  $$
  where $\Delta$ is the modular function\footnote
  {\footcntr}{\eightpoint With the convention that
  $\int_G f(ts) \d t = \Delta(s\inv) \int_G f(t)\d t$, and $\int_G
f(t\inv) \d t = \int_G f(t) \Delta(t\inv) \d t$.}
  of $G$.
  We then see that $\a_s(a)$ is $\a$-integrable.  This clearly implies
that $\Pa$ is $\a$-invariant and, since $\Na$ and $\Ma$ are built in
terms of $\Pa$, it follows that they are also $\a$-invariant.
  We leave the proof of (iv) to the reader.
  \proofend

  \sysstate{Remark}{\rm}
  {\label \RemarkOnAbsValue Regarding
\lcite{\BasicPropertiesOfModule}.(ii) we should remark that Rieffel
\scite{\RieffelIntegrable}{8.9} has found an example of a self-adjoint
$\a$-integrable element $a$ such that $|a|$ is not $\a$-integrable.
This means that $a_+$ and $a_-$ cannot both be $\a$-integrable.  A
short argument involving the fact that $\Pa$ is hereditary shows that
$a$ is not a linear combination of positive $\a$-integrable elements.
Therefore, in Rieffel's example, $\Ma$ is strictly smaller than the
set of all $\a$-integrable elements.}

  The following result, which will be used for our $\Pa$, $\Na$ and
$\Ma$, holds in a much greater generality.  For the sake of
completeness we shall give a detailed proof of it although it has
appeared implicitly in the literature.  See for example
\scite{\Ped}{5.1.2}, \scite{\KR}{7.5.2}, and
\scite{\RieffelIntegrable}{Section 1}.

  \state Proposition
  \label \ConeModuleAlgebra
  Let $\P$ be any hereditary cone in any \cstar-algebra $A$.  Define
  $$
  \N=\{a\in A : a^*a\in \P\}\and \M =\span (\P).
  $$ Then
  \izitem $\N$ is a left ideal in $A$,
  \zitem $\M=\N^*\N$ (linear span of products, no closure),
  \zitem $\M\cap A_+ = \P$,
  \zitem $\M$ is a hereditary *-subalgebra of $A$.

  \proof
  Given $a,b\in \N$ note that
  $$
  (a+b)^*(a+b) \leq (a+b)^*(a+b) + (a-b)^*(a-b) = 2 a^*a + 2 b^*b
  \in \P.
  $$
  Since $\P$ is hereditary we conclude that $(a+b)^*(a+b)\in\P$ and
hence that $a+b\in\N$.  That $\N$ is closed under complex
multiplication is evident, so we see that $\N$ is a linear subspace of
$A$.

  If $c\in A$ and $a\in\N$ we have that
  $
  (ca)^*ca = a^* c^*ca \leq \[c\]^2 a^*a \in \P.
  $
  Since $\P$ is hereditary we see that $ca\in\N$, thus proving that
$\N$ is a left ideal.

  Given $a,b\in\N$ we have by the polarization identity
  $
  a^*b = {1\over 4} \sum_{k=0}^3 i^{-k} (a+i^kb)^*(a+i^kb)
  $
  that $a^*b\in \span(\P)=\M$, and hence that $\N^*\N\subseteq\M$.
Conversely, let $a\in\P$.  Then we have by definition that
$b:=a^{1/2}\in\N$, and hence
  $
  a = b^*b \in \N^*\N.
  $
  This shows that $\P\subseteq\N^*\N$ and hence that
  $\M = \span(\P) \subseteq \N^*\N$.  This proves (ii).

  Speaking of (iii) note that $\M\cap A_+ \supseteq \P$ by definition.
Conversely, let $a\in\M\cap A_+$, and write $a$ as $a = \sum_{i=1}^n
\lambda_i p_i$ with $p_i\in\P$ and $\lambda_i\in\C$.  Since $a =
(a+a^*)/2$ we may assume that the $\lambda_i\in{\bf R}$.  Therefore we
have
  $$
  0 \leq a = \sum_{i=1}^n \lambda_i p_i \leq \sum_{i=1}^n |\lambda_i|
p_i \in\P,
  $$
  from which it follows that $a\in\P$.  This proves (iii).

As for (iv) we have by (ii) that $\M \M = \N^*\N\N^*\N \subseteq \N^*
A\N \subseteq \N^* \N = \M$, where the penultimate step follows by
(i).  This shows that $\M$ is an algebra and the remaining assertions
in (iv) are now evident.
  \proofend

Boosting up the conclusion of \lcite{\BasicPropertiesOfModule}.(iv) we
have:

  \state Proposition
  \label \FirstOcurenceOfRip
  $\Na$ is a
  right pre--Hilbert $\Mult_e(A)$--module
  for the usual multiplication, and inner-product given by
  $$
  \Rip{a}{b} := \sint_G \a_t(a^*b) \d t,
  $$
  for $a,b\in \Na$.

  \proof
  It now suffices to check that
  $$
  \Rip{a}{bm} = \Rip{a}{b}m,
  $$
  for all $a,b\in \Na$ and $m\in \Mult_e(A)$, but this follows easily
by inspection.
  \proofend

  In dealing with Hilbert modules it is always convenient to have a
Hilbert space representation.  The following result is intended to
provide us with one.

  \state Theorem
  \label \ColumnOperators
  For each $a\in\Na$ there exists a bounded linear transformation
  $
  \zeta(a) \colon \H \to L_2(G,\H)
  $
  such that for every $v\in \H$ and $t\in G$,
  $$
  \zeta(a)v \calcat t = \Delta(t)^{-1/2}\a_t\inv (a)v.
  $$
  Moreover the map
  $$
  \zeta : a\in\Na\longmapsto\zeta(a)\in \B(\H, L_2(G,\H))
  $$
  is an isometric representation of the module $\Na$ in the sense that
for any $a,b\in\Na$, and $m\in\Mult_e(A)$,
  \izitem $\[\zeta(a)\] = \[a\]_{\Na}$, where as usual $\[a\]_{\Na} =
\[\Rip{a}{a}\]^{1/2}$,
  \zitem $\zeta(a)m= \zeta(am)$,
  \zitem $\Rip{a}{b} = \zeta(a)^*\zeta(b)$.

  \proof
  Let $a\in\Na$ and $v\in\H$.
  Recall that $A$ is a non-degenerate algebra of operators on $\H$.
Therefore by the Cohen-Hewitt factorization theorem
\scite{\HR}{32.22}, for every $\varepsilon>0$ there exists $b$ in $A$
and $w$ in $\H$ such that $v=bw$, $\[b\]\leq 1$, and
$\[w-v\]<\varepsilon$.  So
  $$
  \[\zeta(a)v \]^2 =
  \int_G \[ \Delta(t)^{-1/2}\a_t\inv (a)v \]^2 \d t =
  \int_G \[ \a_t (a)v \]^2 \d t =
  \int_G \<\a_t (a^*a)bw,v\> \d t \$=
  \<\(\uint_G \a_t (a^*a)b\d t\)w,v\> \leq
  \[\uint_G \a_t (a^*a)b\d t \] \[w\] \[v\] \$\leq
  \[\sint_G \a_t (a^*a) \d t\] \[b\]\[w\] \[v\] \leq
  \[a\]^2_{\Na} (\[v\]+\varepsilon) \[v\].
  $$
  Therefore, since $\varepsilon$ is arbitrary, we have that
  $
  \[\zeta(a)v \] \leq
  \[a\]_{\Na} \[v\] < \infty.
  $
  It follows that $\zeta(a)v$ is indeed in $L_2(G,\H)$, that
$\zeta(a)$ is a bounded map, and that $\[\zeta(a)\] \leq \[a\]_{\Na}$.

  In order to prove (ii) let $m\in\Mult_e(A)$ and note that for all
$t$ in $G$ we have
  $$
  \zeta(a) mv \calcat t =
  \Delta(t)^{-1/2} \a_t\inv (a) mv =
  \Delta(t)^{-1/2}\a_t\inv (am)v =
  \zeta(am)v \calcat t.
  $$
  As for (iii), let $a,b\in\Na$.  Fixing $c\in A$, and $v$ and $w$ in
$\H$, we have
  $$
  \<\zeta(a)^*\zeta(b) cv, w \> =
  \<\zeta(b)cv,\zeta(a)w\> =
  \int_G \Delta(t)^{-1} \<\a_t\inv (b) cv, \a_t\inv (a) w\> \d t \$=
  \int_G \<\a_t(a^*b) cv, w\> \d t =
  \<\(\uint_G\a_t(a^*b) c\d t\)v, w\> =
  \<\Rip ab cv,w\>.
  $$
  Since $A$ is non-degenerate this shows that $\zeta(a)^*\zeta(b) =
\Rip ab$.  It now remains to show the equality in (i) but this follows
easily since
  $$
  \[\zeta(a)\]^2 =
  \[\zeta(a)^*\zeta(a)\] =
  \[\Rip aa\] =
  \[a\]_{\Na}^2.
  \proofend
  $$

Many issues become greatly simplified once we have an isometric
representation of a Hilbert module, such as the one constructed in our
previous result.  For instance:

  \state Corollary
  \label \ModuleProperties
  \izitem The completion of $\Na$, which we henceforth denote by
$\Xa$, can be identified with the closure of $\zeta(\Na)$ within
$\B(\H,L_2(G,\H))$.
  \zitem If $T$ is a bounded linear operator on $L_2(G,\H)$ such that
both
  $T\Xa \subseteq \Xa$ and
  $T^*\Xa \subseteq \Xa$ then the map $S\in\Xa\mapsto TS\in\Xa$ is an
adjointable operator \scite{\Jensen}{1.1.7} on $\Xa$.
  \zitem The algebra $\K(\Xa)$ of generalized compact operators
  \cite{\RieffelInduced}
  on $\Xa$ can be identified with the closed linear span of $\Xa\Xa^*$
within $\B(L_2(G,\H))$.
  \zitem $\Xa$ is a ternary ring of operators \cite{\Zettl} in the
sense that $\Xa \Xa^*\Xa \subseteq \Xa$.

  \proof
  Left to the reader.
  \proofend

  We immediately obtain the following:

  \state Corollary
  Let $D$ be the subalgebra of $\Mult_e(A)$ given by
$D=\overline{\Xa^*\Xa}$ (closed linear span) and let $E$ be the
algebra of operators on $L_2(G,\H)$ given by $E=\overline{\Xa\Xa^*}$.
Then $\Xa$ is an imprimitivity bimodule between $E$ and $D$, which are
therefore {\MR} equivalent \cite{\RieffelInduced}.

We would now like to briefly describe the left regular representation
of the crossed product, mainly to fix our notation.  See \cite{\Ped}
for details. Consider the representation
  $$
  \pi : A \to \B(L_2(G,\H))
  $$
  given by
  $$
  \pi(a)\xi\calcat t = \a_t\inv (a)\xi(t)
  $$
  for all $a\in A$, $\xi\in L_2(G,\H)$ and $t\in G$.  Also, let
$\Lambda$ be the representation of $G$ on $\B(L_2(G,\H))$ given by
$\Lambda=\lambda\* id$, where $\lambda$ is the left regular
representation of $G$ on $L_2(G)$ and we have identified
  $L_2(G,\H) \simeq L_2(G)\* \H$.
  It is well known \scite{\Ped}{7.7.1} that the pair $(\pi,\Lambda)$
is a covariant representation of the \cstar-dynamical system
$(A,G,\a)$ and that
  $$
  \pi\x\Lambda\ :\ \fullcp \longrightarrow \B(L_2(G,\H))
  $$
  is a faithful representation, provided that $G$ is an amenable group
\scite{\Ped}{7.7.5 and 7.7.7}.  In any case $\pi\x\Lambda$ is called
the \stress{regular representation} of $\fullcp$ and its range is the
so called \stress{reduced crossed product} $\rcp$.  We stress that
$\rcp$ is thus a concrete algebra of operators on $L_2(G,\H)$.

Observe that the algebra $E:= \overline{\Xa\Xa^*}$, mentioned above,
is also an algebra of operators on $L_2(G,\H)$.  This raises the
question as to whether there is any relationship between $\rcp$ and
$E$.  This turns out to be the most dramatic question in the present
subject.  We shall have more to say about it in what follows.

For simplicity we let
  $$
  \rho := \pi\x\Lambda.
  $$
  For $f$ in $C_c(G,A)$ (seen as a subalgebra of $\fullcp$),
  $\xi$ in $C_c(G,\H)$ (seen as a subspace of $L_2(G,\H)$),
  and $t$ in $G$, we have \scite{\Ped}{7.7.1},
  $$
  \rho(f)\xi\calcat t =
  \int_G \a_t\inv (f(s))\xi(s\inv t) \d s \$=
  \int_G \a_t\inv (f(ts))\xi(s\inv) \d s =
  \int_G \Delta(s)\inv \a_t\inv (f(ts\inv))\xi(s) \d s.
  $$
  It follows that
  $\rho(f)$ is an ``integral operator'' with ``kernel''
  $$
  k(t,s) = \Delta(s)\inv \a_t\inv (f(ts\inv)),
  \leqno{(\seqnumbering)}
  \label \KerOfRegRep
  $$
  meaning that
  $$
  \rho(f)\xi\calcat t =
  \int_G k(t,s)\xi(s) \d s.
  \leqno{(\seqnumbering)}
  \label \RegRepAsIntOp
  $$
  One easily checks that $k$ satisfies
  $$
  k(tr,sr) = \Delta(r)\inv\a_r\inv (k(t,s)),
  \leqno{(\seqnumbering)}
  \label \LaurentCond
  $$
  for all $t,s,r\in G$.  Later we will investigate integral operators
  again.

  \section {Integrable elements}
  This section is intended to discuss a few technical results about
unconditional integrability to be used below.  Our main reference for
what follows is \cite{\unconditional}.  We begin by quoting a few
consequences of \cite{\unconditional}, in a form suitable for our
purposes.

  \state Theorem
  \label \ResultsFromUnconditional
  Let $f\colon G \to A$ be an unconditionally integrable map.
  \izitem For every $\phi\in L_\infty(G)$ one has that the pointwise
product $\phi f$ is also unconditionally integrable.
  \zitem There exists $M\geq0$ such that
  $$
  \[\uint_G \phi(t)f(t) \d t\] \leq M \[\phi\]
  \for\phi\in L_\infty(G).
  $$
  \zitem For every $\varepsilon>0$ there exists a compact set
$K\subseteq G$ such that for every compact subset $L$ of $G$ with
$K\cap L=\emptyset$,
  $$
  \[\uint_L \phi(t)f(t) \d t\] \leq \varepsilon \[\phi\]
  \for\phi\in L_\infty(G).
  $$

  \proof The first statement is \scite{\unconditional}{2.8}.  The
second follows easily from (i) and \scite{\unconditional}{2.7}.
Finally, (iii) is precisely \scite{\unconditional}{2.9}.
  \proofend

  \state Proposition
  \label \PropositionToDefineOneNorm
  Let $a\in A$ be $\a$-integrable.  Then for every $\phi\in
L_\infty(G)$ the map
  $
  t\mapsto \phi(t)\a_t(a)
  $
  is strictly-unconditionally integrable.  In addition there exists a
constant $M\geq0$ such that
  $$
  \[\sint_G \phi(t)\a_t(a) \d t\] \leq M \[\phi\]
  \for \phi\in L_\infty(G).
  $$

  \proof
  For each $\phi\in L_\infty(G)$ consider the maps
  $$
  L_\phi(b) = \uint_G \phi(t)\a_t(a)b \d t
  \and
  R_\phi(b) = \uint_G \phi(t)b\a_t(a) \d t
  \for b\in A.
  $$
  It is clear that the pair $(L_\phi,R_\phi)$ defines a multiplier of
$A$.  Since multipliers are automatically bounded we have, in
particular, that
  $L_\phi\in \B(A,A)$.  We claim that the set $\{L_\phi\colon
\[\phi\]\leq1\}$ is pointwise bounded.  In fact, fixing $b$ in $A$,
the map $t\mapsto \a_t(a)b$ is unconditionally integrable and hence by
\lcite{\ResultsFromUnconditional}.(ii) there exists a constant $M$
such that
  $$
  \[L_\phi(b)\] =
  \[\uint_G \phi(t)\a_t(a)b \d t\] \leq M\[\phi\] \leq M,
  $$
  provided that $\[\phi\]\leq1.$ By the uniform boundedness principle
there exists a constant $N$ such that $\[L_\phi\]\leq N\[\phi\]$ for
all $\phi\in L_\infty(G)$ and hence,
  $$
  \[\uint_G \phi(t)\a_t(a)b \d t\] =
  \[L_\phi(b)\] \leq
  \[L_\phi\]\[b\] \leq
  N \[\phi\]\[b\],
  $$
  for all $b\in A$ and $\phi\in L_\infty(G)$.  This concludes the
proof.
  \proofend

  \definition
  Given an $\a$-integrable element $a\in A$ we shall denote by
$\[a\]_1$ the smallest constant $M$ for which the inequality in
\lcite{\PropositionToDefineOneNorm} holds.

If $a$ is $\a$-integrable one could attempt to introduce a different
``$L_1$-norm'' of $a$ by setting
  $
  \[a\]_1':=\[\sint_G\a_t(|a|)\d t\].
  $
  However this does not work because $|a|$ may not be $\a$-integrable
as remarked in \lcite{\RemarkOnAbsValue}.  See also Rieffel's
observation following \scite{\RieffelIntegrable}{1.1}

  \state Proposition
  \label \UnifOnCpcts
  Let $f\colon G\to A$ be an unconditionally integrable map and let
$\{\phi_i\}_i\subseteq L_\infty(G)$ be a bounded net converging to
$\phi \in L_\infty(G)$ uniformly over compact subsets of $G$.  Then
  $$
  \lim_i\uint_G \phi_i(t) f(t) \d t =
  \uint_G \phi(t) f(t) \d t
  $$
  in the norm topology of $A$.

  \proof
  Given $\varepsilon>0$ let $M$ be as in
\lcite{\ResultsFromUnconditional}.(ii) and $K$ as in
\lcite{\ResultsFromUnconditional}.(iii).
  Choose an index $i_0$ such that for all $i\geq i_0$ one has that
  $
  \sup_{t\in K} \[\phi_i(t)-\phi(t)\] \leq \varepsilon.
  $
  Therefore for $i\geq i_0$
  $$
  \[\uint_G \phi_i(t) f(t) \d t - \uint_G \phi(t) f(t) \d t\] \$\leq
  \[\int_K (\phi_i(t)-\phi(t)) f(t) \d t\] +
  \[\uint_{G\setminus K} (\phi_i(t)-\phi(t)) f(t) \d t\] \$\leq
  M\[\char K (\phi_i-\phi)\] + \varepsilon \[\phi_i-\phi\] \leq
  M\varepsilon + 2 \varepsilon \sup_i\[\phi_i\].
  $$
  This concludes the proof.
  \proofend

  \state Proposition
  \label \Ineq
  Let $f\in C_c(G,A)$.  Then
  $$
  \(\int_G f(t) \d t\)^* \(\int_G f(t) \d t \)\leq
  |\text{supp}(f)|
  \int_G f(t)^*f(t) \d t,
  $$
  where $|\text{supp}(f)|$ refers to the Haar measure of the support
of $f$.

  \proof
  Let $S$ be the support of $f$ and recall that $A$ is an operator
algebra on the Hilbert space $H$.  For $v\in\H$ we have

  $$
  \<\(\int_S f(t) \d t\)^{\hbox{$*$}}\(\int_S f(t) \d t\) v,v\> =
  \[\int_S f(t)v \d t\]^2 \leq
  \( \int_S \[f(t)v\] \d t \)^2 \$=
  \(\int_S \<f(t)v,f(t)v\>^{1/2} \d t\)^2 =
  \(\int_S \<f(t)^*f(t)v,v\>^{1/2} \d t\)^2 \$\leq
  \(\int_S 1 \d t \)
  \( \int_S \<f(t)^*f(t)v,v\> \d t \) =
  |\text{supp}(f)|
  \<\(\int_S f(t)^*f(t)\d t\)v,v\>,
  $$
  where the penultimate step is H\"older's inequality.  Since $v$ is
arbitrary, the proof is concluded.
  \proofend

  \state Proposition
  \label \IntIsIntegrable
  Let $a$ be an $\a$-integrable element of $A$ and let $g\in C_c(G)$.
Then the element $a'\in A$ defined by
  $
  a' = \int_G g(t)\a_t(a) \d t
  $
  is $\a$-integrable as well.  In addition, if $\phi\in L_\infty(G)$
then
  $$
  \sint_G \phi(t) \a_t(a') \d t =
  \sint_G (\phi\conv g)(t) \a_t(a) \d t.
  $$

  \proof
  For the first assertion we have to show that the net
  $$
  \left\{\int_K b\a_s(a')c \d s \right\}_{K\in\K}
  $$
  converges in the norm of $A$, whenever $b\in A$ and $c=1$, or $b=1$
and $c\in A$.  We have
  $$
  \int_K b\a_s(a')c \d s =
  \int_K b\a_s\( \int_G g(t) \a_t(a) \d t\)c \d s =
  \int_K \int_G g(t) b\a_{st}(a)c \d t \d s \$=
  \int_K \int_G g(s\inv t) b\a_t(a)c \d t \d s =
  \int_G \(\int_K g(s\inv t) \d s\) b\a_t(a)c \d t \$=
  \int_G \(\int_G \char K(s) g(s\inv t) \d s\) b\a_t(a)c \d t =
  \int_G (\char K * g) (t) b\a_t(a)c \d t.
  $$
  The proof will be concluded through an application of
\lcite{\UnifOnCpcts}, once we show that the net
  $\{\char K * g\}_{K\in\K}$ is bounded and converges uniformly over
compacts to the constant function taking the value $\int_G g(s\inv)\d
s$ throughout $G$.  For this purpose let $S$ be the support of $g$.
Given an arbitrary compact subset $L \subseteq G$, set $K_0=LS\inv$.
We claim that if $K\supseteq K_0$, and $t\in L$, then $(\char K *
g)(t) = \int_G g(s\inv)\d s$.  In fact,
  $$
  (\char K * g)(t) =
  \int_G \char K(s) g(s\inv t) \d s =
  \int_G \char K(ts) g(s\inv) \d s =
  \int_{t\inv K} g(s\inv) \d s,
  $$
  but since $t\in L$, we have,
  $$
  tS\inv \subseteq
  LS\inv = K_0 \subseteq K.
  $$
  So $S\inv \subseteq t\inv K$ and hence
  $
  \int_{t\inv K} g(s\inv) \d s =
  \int_{S\inv} g(s\inv) \d s =
  \int_G g(s\inv) \d s.
  $
  This concludes the proof of the fact that $a'$ is $\a$-integrable.
As for the last assertion in the statement, it can now be proved by a
simple change of variable, since we now know that both integrals
converge.
  \proofend

  \state Proposition
  \label \IntegralFormInModule
  Let $a\in\Na$ and $g\in C_c(G)$.  Then
  $ \int_G g(t)\a_t(a) \d t\in \Na$.

  \proof
  By \lcite{\Ineq} we have that
  $$
  \(\int_G g(t)\a_t(a) \d t\)^{\hbox{$*$}} \(\int_G g(t)\a_t(a)\d t\)
\leq
  |\text{supp}(g)| \int_G |g(t)|^2 \a_t(a^*a) \d t.
  $$
  Thus, to prove the statement it suffices to show that
  $$
  a':= \int_G |g(t)|^2 \a_t(a^*a) \d t
  $$
  is $\a$-integrable, since the set of positive $\a$-integrable
elements forms a hereditary cone by
\lcite{\BasicPropertiesOfModule}.(i).  But the integrability of $a'$
follows from \lcite{\IntIsIntegrable}.
  \proofend

  \section {The left structure of $\Xa$}
  \label \LeftStructure
  We remain under the assumption that $(A,G,\a)$ is a \cstar-dynamical
system, where $A$ is a non-degenerate \cstar-algebra of operators on a
Hilbert space $\H$, and $G$ is a locally compact group.  We would now
like to relate the space $\Xa$ (introduced in
\lcite{\ModuleProperties}.(i)) to the reduced crossed product algebra
$\rcp$.

Recall that $(\pi,\Lambda)$,
  defined near the end of section \lcite{\BasicFactsSection},
  is a covariant representation of the \cstar-dynamical system
$(A,G,\a)$ on the Hilbert space $L_2(G,\H)$, and that $\rho$ denotes
the representation of $\fullcp$ given by $\rho = \pi\x\Lambda$.

  \state Lemma
  \label \NewActionOfCrossedProduct If $a\in\Na$ then
  \izitem $\pi(b)\zeta(a)=\zeta(ba)$ for all $b \in A$.
  \zitem $\Lambda_t \zeta(a) = \Delta(t)^{1/2} \zeta(\a_t(a))$ for all
$t \in G$.
  \zitem Let $g\in C_c(G)$ and define
  $
  a'= \int_G g(t) \Delta(t)^{1/2} \a_t(a) \d t.
  $
  Then
  $
  \zeta(a') = \(\int_G g(t) \Lambda_t \d t\) \zeta(a)
  $
  (observe that
  $a'\in\Na$ by \lcite{\IntegralFormInModule}).
  \zitem Let $g$ and $a'$ be as above and take $c\in A$.  Consider the
function $f\in L_1(G,A)$ defined by $f(t)=g(t)c$.
  Then $\rho(f)\zeta(a) = \zeta(ca')$
  (observe that
  $ca'\in\Na$ by \lcite{\ConeModuleAlgebra}.(i)).

  \proof
  Let $v\in\H$ and $s\in G$.  Then
  $$
  \pi(b)\zeta(a) v \calcat s =
  \a_s\inv (b)\(\zeta(a) v \calcat s\) =
  \Delta(s)^{-1/2}\a_s\inv (b)\a_s\inv (a)v \$=
  \Delta(s)^{-1/2}\a_s\inv (ba)v =
  \zeta(ba) v \calcat s,
  $$
  proving (i).  As for (ii) we have
  $$
  \Lambda_t \zeta(a) v \calcat s =
  \zeta(a) v \calcat {t\inv s} =
  \Delta(t\inv s)^{-1/2} \a_{t\inv s}\inv (a) v \$=
  \Delta(t)^{1/2} \Delta(s)^{-1/2} \a_s\inv (\a_t(a)) v =
  \Delta(t)^{1/2} \zeta(\a_t(a)) v \calcat s.
  $$
  With respect to (iii) it is certainly tempting to apply $\zeta$ to
both sides in the definition of $a'$ and use (ii) but this would
require some sort of continuity property for $\zeta$ which we do not
seem to have.  Instead, let $v\in H$ and $\xi\in L_2(G,\H)$.  Then
  $$
  \<\(\int_G g(t) \Lambda_t \d t\)\zeta(a) v,\xi\> =
  \int_G \<\(\int_G g(t) \Lambda_t \d t\)\zeta(a) v \calcat s,
\xi(s)\>
  \d s \$=
  \int_G \< \int_G g(t) \(\zeta(a) v\calcat{t\inv s}\)\d t, \xi(s)\>
\d s \$=
  \int_G \< \int_G g(t)
    \Delta(t\inv s)^{-1/2} \a_{t\inv s}\inv(a)v
    \d t, \xi(s)\> \d s \$=
  \int_G \int_G g(t) \Delta(t)^{1/2}\Delta(s)^{-1/2}
  \<\a_{s\inv t}(a)v, \xi(s)\> \d t\d s.
  $$
  On the other hand
  $$
  \< \zeta(a') v,\xi\> =
  \int_G \<\zeta(a') v\calcat s,\xi(s)\> \d s =
  \int_G \<\Delta(s)^{-1/2}\a_s\inv (a')v,\xi(s)\> \d s \$=
  \int_G \<\Delta(s)^{-1/2}\a_s\inv \( \int_G g(t) \Delta(t)^{1/2}
\a_t(a) \d t \)v,\xi(s)\> \d s \$=
  \int_G \int_G
  \Delta(s)^{-1/2}g(t) \Delta(t)^{1/2}\<\a_{s\inv t}(a)v,\xi(s)\> \d
t\d s.
  $$
  This proves (iii).  As for (iv) we have
  $$
  \rho(f) \zeta(a) =
  \(\int_G \pi(f(t))\Lambda_t \d t\) \zeta(a) =
  \pi(c)\(\int_G g(t) \Lambda_t \d t\) \zeta(a) =
  \pi(c)\zeta(a') = \zeta(ca'),
  $$
  proving (iv).
  \proofend

  \state Corollary
  \label \LeftCrossedProductModule
  Let $\N$ be a subset of $\Na$ such that
  \izitem for every $g\in C_c(G)$ and $a\in\N$ one has that $\int_G
g(t)\a_t(a)\d t \in\N$, and
  \zitem there exists a dense subset ${\cal D}$ of $A$ such that
${\cal D}\N\subseteq\N$.
  \medskip\noindent
  Then $\(\rcp\) \X \subseteq \X$, where $\X$ is the closure of
$\zeta(\N)$ in $\Xa$.

  \proof
  Let $f(t) = g(t)c$, where $g\in C_c(G)$ and $c\in \cal D$.  Then, by
\lcite{\NewActionOfCrossedProduct}.(iv) we have that
  $\rho(f)\zeta(\N) \subseteq \zeta(\N)$, and hence that
  $\rho(f)\X \subseteq \X$.  Since the linear span of the set of all
$\rho(f)$'s of the above form is dense in $\rcp$, the conclusion
follows.
  \proofend

  \state Corollary
  One has that
  $\(\rcp\) \Xa \subseteq \Xa$

  \proof
  Follows at once from \lcite{\LeftCrossedProductModule} once we note
that $\Na$ satisfies the required hypothesis by
\lcite{\IntegralFormInModule} and
  \lcite{\ConeModuleAlgebra}.(i).
  \proofend

We therefore see that $\Xa$ is a left $\(\rcp\)$--module under the
composition of operators.  Moreover, in view of
\lcite{\ModuleProperties}.(ii), all operators in $\rcp$ act as
adjointable operators on $\Xa$.

If we let $E= \overline{\Xa\Xa^*}$, as before, we have that
  $\(\rcp\) E \subseteq E$ and hence that
  $\(\rcp\)\cap E$ is an ideal in $\rcp$.  It would seem natural to
conjecture that this is the ideal one would like to show is {\MR}
equivalent to the (still not yet defined) {\GFPA}, as in
\scite{\RieffelIntegrable}{Section 6}.
  A particularly intriguing question is:

  \sysstate{Question}{\rm}
  {Is there a linear subspace $\X$ of $\Xa$ such that
$\overline{\X\X^*} =\(\rcp\)\cap E$?}

A related question is whether or not $\Xa$ is also a Hilbert module
over $\rcp$.  Clearly this would be the case if we knew that $\rcp$
contains the range of the left inner-product on $\Xa$, given by
  $$
  \Lip TS = TS^*
  \for T,S\in \Xa.
  $$
  In particular, one could ask:

  \sysstate{Question}{\rm}
  {Given a pair of elements $a,b\in\Na$, how could one determine if
$\zeta(a)\zeta(b)^*$ belongs to $\rcp$?}

This question, to which we will give a satisfactory answer in the
abelian group case, resides in the heart of the matter and will
dominate our attention in the remaining sections of this work.  In
order to explore it further we need a deeper understanding of integral
operators.

  \section {Integral and Laurent operators}
  Inspired by \lcite{\RegRepAsIntOp} we would now like to establish a
precise notion of integral operators in our context.  For this purpose
let $k$ be any continuous $\B(\H)$-valued function on $G\x G$ and
suppose that for all $\xi\in C_c(G,\H)$ the function
  $\eta : G \to \H$, given by
  $$
  \eta(t) = \int_G k(t,s)\xi(s)\d s
  \for t\in G,
  $$
  is in $L_2(G,\H)$.  Suppose further that there exists a constant
$M\geq 0$ such that $\[\eta\]_2 \leq M \[\xi\]_2$ for every $\xi$ as
above.  This implies the existence of a bounded operator $T$ on
$L_2(G,\H)$ satisfying
  $$
  T\xi\calcat t =
  \int_G k(t,s)\xi(s) \d s
  \for \xi\in C_c(G,\H)\for t\in G.
  \leqno{(\seqnumbering)}
  \label \IntegralOp
  $$

  \definition
  By an \stress{integral operator} on $L_2(G,\H)$ we shall mean any
bounded operator $T$ that satisfies \lcite{\IntegralOp} for some
continuous function $k\colon G\x G \to \B(\H)$.  In this case $k$ will
be called the \stress{kernel} of $T$.

  If a general theory of integral operators is desired one should
probably relax the requirement that the integral kernel $k$ be
continuous.  For the applications we have in mind, however, the
definition given is enough, apart from the fact that it greatly
simplifies the study of our operators.

  Whereas it is highly unlikely that a necessary and sufficient
condition on $k$ will ever be found for $T$ to be bounded, we will be
able to use the concept of integral operators quite profitably.  Of
course, in the applications, the boundedness of $T$ must be derived
{}from other sources.  This is akin to studying an operator on a Hilbert
space via its matrix with respect to an orthonormal basis although no
one really knows how to characterize boundedness in terms of matrices.

  Observe that \lcite{\RegRepAsIntOp} shows that $\rho(f)$ is an
integral operator for any $f$ in $C_c(G,A)$.

  Inspired by \lcite{\LaurentCond}, and in analogy with the case of
the trivial action of ${\bf Z}$ on ${\bf C}$, we make the following:

  \definition
  A \stress{Laurent operator} is an integral operator $T$ on
$L_2(G,\H)$ whose kernel $k$ takes values in $A\subseteq\B(\H)$ and
satisfies \lcite{\LaurentCond}, namely that
  $k(tr,sr) = \Delta(r)\inv \a_r\inv (k(t,s))$
  for all $t,s,r\in G$.

  Given a Laurent operator with kernel $k$, define
  $$
  f(r):=\a_r(k(r,e)) \for r\in G.
  $$
  Taking into account the right hand side of \lcite{\KerOfRegRep}, let
us compute
  $$
  \Delta(s)\inv \a_t\inv (f(ts\inv)) =
  \Delta(s)\inv \a_t\inv (\a_{ts\inv} (k(ts\inv,e))) =
  \Delta(s)\inv \a_s\inv (k(ts\inv,e)) =
  k(t,s).
  $$
  Therefore $k$ satisfies \lcite{\KerOfRegRep} with respect to $f$.
The reader should however be warned that, in general, the function $f$
defined above for a Laurent operator need not be in $C_c(G,A)$, or
even in $L_1(G,A)$.  This is related to the remark at the end of
section 6 in \cite{\RieffelIntegrable}.

  \definition
  Given a Laurent operator $T$ with kernel $k$ we say that the
(continuous) function
  $f\colon G \to A$ given by $f(r)=\a_r(k(r,e))$ is the
\stress{symbol} of $T$.

  It is then obvious that for any $f\in C_c(G,A)$ the symbol of
$\rho(f)$ is $f$.  We also note that:

  \state Proposition
  \label \NecCondOnSymbol
  Let $T$ be a Laurent operator with symbol $f$.  Suppose that $f$ is
in $C_c(G,A)$.  Then $T$ belongs to $\rcp$.

  \proof
  By \lcite{\IntegralOp} it is clear that two Laurent operators with
the same symbol, and hence also the same kernel, coincide.  Since both
$T$ and $\rho(f)$ are Laurent operators with symbol $f$, we must have
$T=\rho(f)$.
  \proofend

  One of the reasons we are interested in Laurent operators is as
follows:

  \state Proposition
  \label \AtLeastLaurent
  Given $a,b\in\Na$ one has that $\zeta(a)\zeta(b)^*$ is a Laurent
operator with symbol
  $$
  f(r) = \Delta(r)^{-1/2}a \a_r(b^*).
  $$

  \proof
  Let $\xi\in C_c(G,\H)$ and $v\in\H$.  Then
  $$
  \<\zeta(b)^*\xi,v\> =
  \<\xi,\zeta(b)v\> =
  \int_G \<\xi(s),\Delta(s)^{-1/2} \a_{s}\inv(b) v\> \d s \$=
  \int_G \<\Delta(s)^{-1/2} \a_{s}\inv(b^*)\xi(s),v\> \d s.
  $$
  So we see that $\zeta(b)^*\xi$ is given by
  $$
  \zeta(b)^*\xi = \int_G \Delta(s)^{-1/2} \a_{s}\inv(b^*)\xi(s) \d s,
  $$

  Given $t\in G$ we then have
  $$
  \zeta(a)\zeta(b)^*\xi \calcat t =
  \Delta(t)^{-1/2} \a_t\inv (a)
  \zeta(b)^*\xi =
  \Delta(t)^{-1/2} \a_t\inv (a)
  \int_G \Delta(s)^{-1/2} \a_{s}\inv(b^*)\xi(s) \d s \$=
  \int_G
    \Delta(ts)^{-1/2}\a_{t}\inv(a)\a_{s}\inv(b^*)
    \xi(s)\d s,
  $$
  from where we conclude that $\zeta(a)\zeta(b)^*$ is the integral
operator with kernel
  $$
  k(t,s) =
  \Delta(ts)^{-1/2}\a_{t}\inv(a)\a_{s}\inv(b^*).
  $$
  It is easy to see that $k$ satisfies \lcite{\LaurentCond} and hence
that
  $\zeta(a)\zeta(b)^*$ is a Laurent operator.  Its symbol is given by
  $$
  f(r) =
  \a_r(k(r,e)) =
  \Delta(r)^{-1/2}a \a_r(b^*).
  \proofend
  $$

  Referring to the question posed at the end of section
\lcite{\LeftStructure}, namely whether $\zeta(a)\zeta(b)^*\in\rcp$ for
$a,b\in\Na$, we may give an affirmative answer in a very simple case:

  \state Proposition
  \label\Simplecase
  Suppose $a,b\in\Na$ are such that the map $r\mapsto a\a_r(b^*)$ is
compactly supported.  Then $\zeta(a)\zeta(b)^*\in\rcp$.

  \proof
  This is an immediate consequence of \lcite{\AtLeastLaurent} and
\lcite{\NecCondOnSymbol}.
  \proofend

  \section {Abelian groups and the Fourier transform}
  In the general case there is not much more we can say about the
question of whether $\zeta(a)\zeta(b)^*$ belongs to $\rcp$ for $a,b\in
\Na$, as mentioned at the end of Section \lcite{\LeftStructure}.  This
is the main obstacle to defining the {\GFPA} and proving it to be
{\MR} equivalent to an ideal in $\rcp$.

  We shall therefore restrict our attention to the special case of
abelian groups and hence we assume from now on that $G$ is abelian.
  In particular $G$ is amenable and hence $\rho$ establishes an
isomorphism between $\fullcp$ and $\rcp$ \scite{\Ped}{7.7.7}.
  We will therefore identify these algebras without further notice,
remarking, however, that we will be much more interested in the
concrete algebra $\rcp$ of operators on $L_2(G,\H)$ rather than in the
abstract \cstar-algebra $\fullcp$.

  Let $\Gdual$ be the Pontryagin dual of $G$.  Given $x$ in $\Gdual$
and $t$ in $G$ we will denote the value of the character $x$ on $t$ by
$\tdu {x}{t}$.

  \definition
  Let $a\in A$ be $\a$-integrable and let $x\in \Gdual$.  The
\stress{Fourier coefficient}\footnote{\footcntr}
  {\eightpoint
  \global\edef\Convention{\number \halffootno}
  When working with abelian groups one often has to make choices, as
is the case of the complex conjugation in the definition of the
Fourier transform: one could very well develop all of classical
Harmonic Analysis defining the Fourier transform of a complex function
$f$ by $\^f(x) = \int_G \du xt f(t) \d t$, as opposed to the more
usual $\^f(x) = \int_G \idu xt f(t) \d t$.
  Likewise, given an action of $G$ on a space $X$ one usually induces
an action on the algebra of functions on $X$ by the formula
  $\a_t(f)(x) = f(t\inv x)$, but when $G$ is commutative one has the
option of dropping the inversion of $t$.  This extra freedom has a
price: although there are no right or wrong choices, some of them
often cause an excessive number of inverses and complex conjugations,
which one feels should not be there.
  Often an attempt to back up and change conventions reveals only too
late that the undesirable inverses and conjugations pop up in greater
numbers further on.
  \hfill\break\indent
  The convention adopted here for the Fourier transform (without the
complex conjugation) takes into account that the action of a compact
abelian group $G$ on itself by left multiplication, once induced to
$C(G)$ via the formula
  $\a_t(f)(x) = f(t\inv x)$, ought to satisfy $E_x(f)=\delta_{x,y}f$,
when $f(\cdot) = \du y\cdot$, hence agreeing with the prevailing
convention for the Fourier transform (with the complex conjugation).
This in turn seems to indicate the need to include the complex
conjugation in the definition of
  $\Mult_x(A)$ above.  These choices seem rather reasonable but they
have the somewhat unpleasant consequence of forcing us to part with
tradition with respect to the dual action (see e.g
\scite{\Ped}{7.8.3}) to be defined in section
\fcite{10}{SectionOnDualAction}.
  \hfill\break\indent
  In fact there are stronger reasons for adopting our conventions:
among these we would like to mention that Proposition
\fcite{11.12}{PropositionOfDiagram} would have to suffer some rather
unnatural modifications in order to survive under seemingly natural
alternative conventions.  The same goes for Proposition
\fcite{6.5}{FourierTransfOfConvolution.}}
  $\F ax$ of $a$ is the element of $\Mult(A)$ given by
  $$
  \F ax = \sint_G \du{x}{t}\a_t(a)\d t.
  $$

It is easy to prove \scite{\unconditional}{6.4} that $\F ax$ belongs
to the \stress{$x$-spectral subspace} of $\Mult(A)$ defined by
  $$
  \Mult_x(A) = \{ m \in \Mult(A) : \a_t(m) = \idu xt m,\ t\in G\}.
  $$
  In particular $\Mult_e(A)$ is consistent with our previous notation
for the set of fixed points for $\a$ within $\Mult(A)$.

  As an immediate consequence of \lcite{\PropositionToDefineOneNorm}
we have:

  \state Proposition
  \label \NormOfFourierVsOneNorm
  Let $a$ be $\a$-integrable.  Then for every $x\in\Gdual$ one has
that $\[\F ax\]\leq \[a\]_1$.

  As in the classical case we have:

  \state Proposition
  \label \FourierIsStrictlyContinuous
  {\rm (See also \scite{\unconditional}{6.3}).}
  For each $\a$-integrable element $a\in A$ the \stress{Fourier
transform}
  $$
  x \in \Gdual \mapsto \F ax \in \Mult(A)
  $$
  is uniformly continuous in the strict topology of $\Mult(A)$.  That
is, given $b,c\in\Mult(A)$ such that either $b\in A$ and $c=1$, or
$b=1$ and $c\in A$, one has
  $$
  \lim_{z\to e} \sup_{x\in\Gdual} \[b\F a{zx} c - b\F a{x} c\] =0.
  $$

  \proof
  Suppose by contradiction that this is not so.  Then there exists
$\varepsilon>0$ and nets $\{z_i\}_i$ and $\{x_i\}_i$ in $\Gdual$ such
that $z_i\to 0$ and
  $
  \[b\F a{z_ix_i} c - b\F a{x_i} c\]\geq \varepsilon.
  $
  Observe that
  $$
  b\F a{z_ix_i} c - b\F a{x_i} c =
  \uint_G (\du {z_ix_i}t - \du {x_i}t)b\a_t(a)c \d t \$=
  \uint_G (\du {z_i}t - 1)\du {x_i}t b\a_t(a)c \d t =
  \uint_G \varphi_i(t) b\a_t(a)c \d t,
  $$
  where $\varphi_i(t)=(\du {z_i}t - 1)\du {x_i}t$.
  By definition of the topology on $\Gdual$ we have that $\tdu{z_i}t
\to 1$ uniformly over compact sets and hence $\varphi_i(t)\to0$, also
uniformly over compacts.
  Using \lcite{\UnifOnCpcts} we thus obtain
  $$
  \int_G \varphi_i(t) b\a_t(a)c \d t \to 0,
  $$
  therefore arriving at a contradiction.
  \proofend

For future use we now collect some important properties of the Fourier
transform.

  \state Proposition
  \label \FormulasForFourier
  Let $a,b\in A$ be $\a$-integrable, let $x,y\in\Gdual$, and let $m\in
\Mult_y(A)$.  Then
  \izitem $\F ax ^* = \F {a^*}{x\inv}$,
  \zitem $ma$ is $\a$-integrable and $m\F ax = \F{ma}{yx}$,
  \zitem $am$ is $\a$-integrable and $\F ax m = \F{am}{xy}$,
  \zitem $\F ax \F by = \F{a \F by}{xy} =\F{\F ax b}{xy}.$

  \proof
  Left to the reader.
  \proofend

Our next result is related to the classical result according to which
the Fourier transform of a convolution is the pointwise product of the
corresponding transforms.  As usual we denote by $\^g$ the Fourier
transform of a function $g\in L_1(G)$, i.e,
  $$
  \^g(x) = \int_G \idu xt g(t) \d t \for x\in\Gdual.
  $$

  \state Proposition
  \label \FourierTransfOfConvolution
  Given an $\a$-integrable element $a\in A$ and $g\in C_c(G)$ let
  $$
  a'=\int_G g(t) \a_t(a) \d t.
  $$
  Then for every $x$ in $\Gdual$ one has that $\F {a'}x = \^g(x) \F
ax$.

  \proof
  Recall initially that $a'$ is $\a$-integrable by
\lcite{\IntIsIntegrable}.  Also from \lcite{\IntIsIntegrable} we have
  $$
  \F {a'}x =
  \sint_G \du xt \a_t(a') \d t =
  \sint_G (\phi\conv g)(t) \a_t(a) \d t,
  $$
  where $\phi(t)= \du xt$.  On the other hand
  $$
  (\phi\conv g)(t) =
  \int_G \phi(ts\inv) g(s) \d s =
  \int_G \du xt \idu xs g(s) \d s =
  \du xt\ \^g(x).
  $$
  Therefore
  $$
  \F {a'}x =
  \sint_G \du xt\ \^g(x) \a_t(a) \d t =
  \^g(x) \F ax.
  $$
  \proofend

The following is the Fourier inversion Theorem for our context.

  \state Proposition
  \label \InvFourierTransform
  Let $a$ be an $\a$-integrable element of $A$ whose Fourier transform
is absolutely integrable, that is, such that $\int_\Gdual \[\F ax\] \d
x <\infty$.  Then for every $t$ in $G$ one has that
  $$
  \int_{\Gdual} \idu xt \F ax \d x = \a_t(a).
  $$

  \proof
  Initially observe that the map $x \mapsto \F ax$ is continuous in
the strict topology of $\Mult(A)$ by
\lcite{\FourierIsStrictlyContinuous}, and hence the integral in the
statement is well defined as a strict integral.  Therefore, in order
to prove the statement, it is enough to show that for every $b\in A$
and vectors $\xi,\eta\in\H$ (recall that $A$ is represented as a
non-degenerate algebra of operators on $\H$) one has that
  $$
  \int_{\Gdual} \idu xt \<\F ax b\xi,\eta\>\d x =
\<\a_t(a)b\xi,\eta\>.
  $$ Consider the function
  $$
  \psi: t\in G \mapsto \<\a_t(a) b\xi, \eta \> \in \C.
  $$
  Since $a$ is $\a$-integrable we have that $\psi$ is in $L_1(G)$.
The inverse Fourier transform of $\psi$ is clearly given by
  $
  \check\psi(x) = \<\F ax b\xi,\eta\>.
  $
  By hypothesis $\check\psi$ is integrable and hence the result
follows from the classical Fourier inversion Theorem
\scite{\HR}{31.44.(c)}.
  \proofend

  \section {The dual unitary group}
  For a given $x$ in $\Gdual$ consider the unitary operator $\V_x$ on
$L_2(G,\H)$ defined by
  $$
  \V_x\xi\calcat t =
  \idu {x}{t}\xi(t)
  \for \xi\in L_2(G,\H), \quad t\in G.
  $$
  It is well known that the correspondence $x\mapsto \V_x$ is a
strongly continuous unitary representation of $\Gdual$ on $L_2(G,\H)$
which we shall call the \stress{dual unitary group}.

  Let $T$ be an integral operator with kernel $k$.  Then for every
$\xi$ in $C_c(G,\H)$ and $t\in G$ we have
  $$
  \V_x T \V_x\inv \xi \calcat t =
  \idu xt \int_G k(t,s) \idu{x\inv}{s} \xi(s) \d s=
  \int_G \idu{x}{ts\inv} k(t,s) \xi(s) \d s,
  $$
  and hence we see that $\V_x T \V_x\inv$ is an integral operator
whose kernel is given by $k_x(t,s) = \idu{x}{ts\inv} k(t,s)$.  It is
also evident that if $T$ is a Laurent operator then so is $\V_x T
\V_x\inv$.  In this case let $f$ be the symbol of $T$.  Then the
symbol $f_x$ of $\V_x T \V_x\inv$ is given by
  $$
  f_x(r) = \a_r(k_x(r,e)) =
  \idu{x}{r} \a_r(k(r,e)) =
  \idu{x}{r} f(r)
  \for r\in G.
  $$

In particular, given $f\in C_c(G,A)$ we saw that $\rho(f)$ is the
Laurent operator with symbol $f$ and hence
  $\V_x \rho(f) \V_x\inv$
  is the Laurent operator with symbol
  $$
  \ah_x(f)(r) := \idu{x}{r}f(r)
  \for r\in G.
  $$
  In other words
  $$
  \V_x \rho(f) \V_x\inv = \rho(\ah_x(f)).
  $$
  Observe that $\ah$ is, up to a sign convention, the dual action of
$\Gdual$ on $\fullcp$, as defined in \scite{\Ped}{7.8.3}, and hence we
see that the pair $(\rho,\V)$ is a covariant representation of the
dual \cstar-dynamical system $(\fullcp, \Gdual, \ah)$.

Again using \scite{\Ped}{7.8.3} we have that the dual action is
strongly continuous and hence for each operator $T$ in $\fullcp$ the
map
  $$
  x \in \Gdual \mapsto \V_x T \V_x\inv \in \B(L_2(G,\H))
  $$
  is norm continuous.

  \definition
  \label \DefinitionOfVContinuous
  We say that a bounded operator $T \in \B(L_2(G,\H))$ is
\stress{continuous with respect to $V$}, or just
\stress{$V$-continuous}, if the map given by
  $x \in \Gdual \mapsto \V_x T \V_x\inv \in \B(L_2(G,\H))$
  is continuous in norm.

  So, for any $f$ in $C_c(G,A)$, the operator $\rho(f)$ is a
$\V$-continuous Laurent operator.  We would now like to prove a
converse of this statement, in which we will use the following:

  \state Lemma
  \label \TauberianNet
  There exists a net $\{g_i\}_{i\in\Lambda}\subseteq L_1(G)$ such that
$\^g_i$ has compact support and $\[g_i\]_1=1$ for all $i\in\Lambda$,
and also such that $g_i\conv f$ converges uniformly to $f$ for any
bounded uniformly continuous function $f$ from $G$ to any Banach space
$X$.

  \proof
  Let $\I(G)$ be the ideal in $L_1(G)$ (under convolution) formed by
all $g\in L_1(G)$ such that $\^g$ has compact support.  Since $L_1(G)$
satisfies Ditkin's condition \scite{\HR}{39.29} we have that $\I(G)$
is dense in $L_1(G)$.

  The usual argument of taking functions supported in smaller and
smaller neighborhoods of the unit of $G$ provides a net as required,
except that $\^g_i$ may not have compact support.  For each
  $(i,n)\in \Lambda\x{\bf N}$ choose $g_{i,n}\in\I(G)$ such that
$\[g_i-g_{i,n}\]\leq 1/n$.  One may now easily prove that the net
$\{g_{i,n}\}_{(i,n)\in\Lambda\x{\bf N}}$
  satisfies the required properties.
  \proofend

This brings us to the main result of this section:

  \state Theorem
  \label \MembershipCriteria
  Let $T$ be a Laurent operator.  Then $T$ belongs to $\fullcp$ if and
only if $T$ is continuous with respect to the dual group $V$.

  \proof
  The ``only if'' part has already been verified in the discussion
before \lcite{\DefinitionOfVContinuous} so let us deal with the
converse and hence we suppose that $T$ is a $\V$-continuous Laurent
operator.  Therefore the expression
  $$
  \tau(x) := \V_x T \V_x\inv
  $$
  defines a norm continuous $\B(L_2(G,\H))$-valued function on
$\Gdual$.  For $x,y\in\Gdual$ note that
  $$
  \[\tau(x)-\tau(y)\] =
  \[\V_x T \V_x\inv - \V_y T \V_y\inv \] =
  \[T- \V_x\inv \V_y T \V_y\inv \V_x\] =
  \[ \tau(e) - \tau(x\inv y)\].
  $$
  Therefore the continuity of $\tau$ at the identity group element $e$
implies that $\tau$ is in fact uniformly continuous on $\Gdual$.
  Reversing the roles of $G$ and $\Gdual$ in \lcite{\TauberianNet},
let $\{g_i\}_i$ be a net in $L_1(\Gdual)$ such that $\^g_i$ has
compact support in $G$ and such that
  $
  \tau_i :=
  g_i\conv \tau
  $
  converges uniformly to $\tau$.  Note that
  $$
  \tau_i(x) =
  \int_\Gdual g_i(y) \tau(y\inv x) \d y =
  \int_\Gdual g_i(y) \V_y\inv \V_x T \V_x\inv \V_y \d y =
  \V_x \( \int_\Gdual g_i(y) \V_y\inv T \V_y \d y\)\V_x\inv.
  $$
  Let us further investigate the operator
  $$
  T_i :=
  \tau_i(e)=
  \int_\Gdual g_i(y) \V_y\inv T \V_y \d y
  $$
  appearing above.  For this purpose let $\xi, \eta \in C_c(G,\H)$ and
note that
  $$
  \<T_i \xi, \eta \> =
  \int_\Gdual g_i(y) \<T \V_y\xi, \V_y \eta \> \d y =
  \int_\Gdual g_i(y) \int_G \int_G
  \<k(t,s) \(\V_y\xi\calcat s\), \V_y \eta \calcat t \>
  \d s\d t\d y \$=
  \int_\Gdual \int_G \int_G
  g_i(y) \idu{y}{t\inv s} \<k(t,s) \xi(s), \eta(t) \>
  \d s\d t\d y \$=
  \int_G \int_G
  \^{g_i}(t\inv s) \<k(t,s) \xi(s), \eta(t) \>
  \d s\d t.
  $$
  It follows that $T_i$ is an integral operator with kernel
  $k_i(t,s) := \^{g_i}(t\inv s) k(t,s)$.  It is evident that $k_i$
satisfies \lcite{\LaurentCond} so that $T_i$ is a Laurent operator.
The symbol $f_i$ of $T_i$ may be computed in terms of the symbol $f$
of $T$ as follows:
  $$
  f_i(r) = \a_r(k_i(r,e)) =
  \^{g_i}(r\inv) \a_r(k(r,e)) =
  \^{g_i}(r\inv) f(r).
  $$

  Since $f$ is continuous and $\^g_i\in C_c(G)$ we have that $f_i$ is
in $C_c(G,A)$.  So by \lcite{\NecCondOnSymbol}, $T_i$ belongs to
$\fullcp$.  Finally, given that $\tau_i$ converges uniformly to $\tau$
we have that
  $$
  T = \tau(e) = \lim_i \tau_i(e) = \lim_i T_i,
  $$
  and hence $T\in \fullcp$ as well.
  \proofend

Given $a,b\in \Na$ we have seen in \lcite{\AtLeastLaurent} that
$\zeta(a)\zeta(b)^*$ is a Laurent operator.  We may therefore apply
the above result to determine whether or not this operator belongs to
$\fullcp$.  Before that we need the following:

  \state Lemma
  \label \FourierViaV
  Let $a,b\in\Na$.  Then for all $x\in\Gdual$ we have
  $$
  \zeta(a)^*V_x\zeta(b) = \F {a^*b}x.
  $$

  \proof
  For $v,w\in\H$ and $c\in A$ we have
  $$
  \<\zeta(a)^*V_x\zeta(b)cv,w\> =
  \<V_x\zeta(b)cv,\zeta(a)w\> =
  \int_G \<\idu xt \a_{t\inv}(b)cv, \a_{t\inv}(a)w\> \d t\$=
  \int_G \<\du xt \a_{t}(a^*b)cv,w\> \d t =
  \< \F{a^*b}x cv,w\>.
  $$
  Observe that, since $G$ is abelian, its modular function $\Delta$ is
identically 1 and hence was omitted in the calculation above.
  Since $A$ is non-degenerate the proof is complete.
  \proofend

  The following result is intended to answer the question posed at the
end of section \lcite{\LeftStructure} for the abelian group case,
giving a necessary and sufficient condition for $\zeta(a)\zeta(b)^*$
to belong to $\fullcp$, for a given pair $a,b\in\Na$.

  \state Theorem
  \label \IncredibleFormula
  Let $a,b\in \Na$ and denote $p:=a^*a$ and $q:=b^*b$.  Then the
following are equivalent:
  \izitem $\zeta(a)\zeta(b)^*$ belongs to $\fullcp$,
  \zitem
  $
  \[\F p{xz} \F qy - \F px \F q{zy}\]
  $
  converges to zero,
  uniformly in $x$ and $y$, as $z\to e$,
  \zitem One has that
  $$
  \lim_{z\to e} \[\F pz \E q - \E p \F qz\] = 0
  $$
  and
  $$
  \lim_{z\to e} \[\F pz \F q{z\inv} - \E p \E q\] = 0.
  $$

  \proof
  \begingroup\def\ca{\zeta(a)}\def\cb{\zeta(b)}
  Observe that, by \lcite{\FourierViaV}, for all $x$, $y$, and $z$ in
$\Gdual$ we have
  $$
  \[\F p{xz} \F qy - \F px \F q{zy}\] =
  \[\ca^* V_{xz}\ca \cb^* V_y\cb - \ca^*V_x\ca \cb^* V_{zy}\cb\]
\$\leq
  \[\ca^* V_x\] \[ V_z\ca \cb^* - \ca \cb^* V_z\] \[ V_y\cb \] \$=
  \[\ca\] \[ V_z\ca \cb^* V_z\inv - \ca \cb^*\] \[ \cb \].
  $$
  Suppose that $\zeta(a)\zeta(b)^*\in \fullcp$.  Then by
\lcite{\MembershipCriteria} we have that $\zeta(a)\zeta(b)^*$ is a
$V$-continuous operator and hence (ii) holds.

  That (ii)$\Rightarrow$(iii) follows by taking $x=y=e$ on the one
hand, and $x=e$ and $y=z\inv$ on the other.
  In order to verify (iii)$\Rightarrow$(i) we claim that
  $$
  \lim_{x\to e} \[V_x \zeta(a)\zeta(b)^* V_x^* - \zeta(a)\zeta(b)^*\]
= 0.
  $$
  Observe that the \cstar-identity $\[a^*a\]=\[a\]^2$ implies that
  $\[aa^*a\]=\[a\]^3$ which is the form we choose to evaluate the norm
in the limit above.  We have
  $$
  \Big(V_x\ca\cb^*V_x^*-\ca\cb^*\Big)
  \Big(V_x\ca\cb^*V_x^*-\ca\cb^*\Big)^*
  \Big(V_x\ca\cb^*V_x^*-\ca\cb^*\Big)\$=
  V_x\ca\Big(\cb^*\cb\ca^*\ca-\cb^*V_x^*\cb\ca^*V_x\ca\Big)\cb^*V_x^*\$+
  \ca\Big(\cb^*\cb\ca^*V_x\ca-\cb^*V_x\cb\ca^*\ca\Big)\cb^*V_x^*\$+
  V_x\ca\Big(\cb^*V_x^*\cb\ca^*\ca-\cb^*\cb\ca^*V_x^*\ca\Big)\cb^*\$+
  \ca\Big(\cb^*\cb\ca^*\ca-\cb^*V_x\cb\ca^*V_x^*\ca\Big)\cb^*\$=
  V_x\ca\Big(\E{q}\E{p}-\F{q}{x\inv}\F{p}{x}\Big)\cb^*V_x^*+
  \ca\Big(\E{q}\F{p}{x}-\F{q}{x}\E{p}\Big)\cb^*V_x^*\$+
  V_x\ca\Big(\F{q}{x\inv}\E{p}-\E{q}\F{p}{x\inv}\Big)\cb^*+
  \ca\Big(\E{q}\E{p}-\F{q}{x}\F{p}{x\inv}\Big)\cb^*,
  $$
  Taking adjoints and using \lcite{\FormulasForFourier}.(i) we
conclude from (iii) that the above tends to zero as $z\to e$.  This
proves our claim and hence that the function
  $$
  \tau \ :\ x\in\Gdual\longmapsto \V_x \zeta(a)\zeta(b)^* \V_x^* \in
\B(L_2(G,\H))
  $$
  is continuous at $x=e$.  As in the proof of
\lcite{\MembershipCriteria} this implies that $\tau$ is uniformly
continuous on $\Gdual$ and hence that $\zeta(a)\zeta(b)^*$ is
$V$-continuous.  That $\zeta(a)\zeta(b)^* \in \fullcp$ then follows
{}from \lcite{\AtLeastLaurent} and \lcite{\MembershipCriteria}.
  \endgroup
  \proofend

  \section {Relative continuity}
  As already mentioned, one of the crucial aspects of this subject is
the question of whether or not $\zeta(a)\zeta(b)^*\in \fullcp$, for
given $a$ and $b$ in $\Na$.  Theorem \lcite{\IncredibleFormula}, which
gives necessary and sufficient conditions for this to happen, will
therefore acquire a special relevance to us.  We will see that a deep
understanding of this issue is the basis for the study of a certain
{\MR} equivalence as well as the spectral theory which we plan to
develop.

  \definition
  Let $a,b\in A$ be $\a$-integrable elements.  We say that the pair
$(a,b)$ is \stress{relatively continuous}, and denote it by $a\rc b$,
if
  $
  \[\F a{xz} \F by - \F ax \F b{zy}\] \to 0
  $
  uniformly in $x$ and $y$, as $z\to e$.

Exercising the new terminology we have the following immediate
consequence of \lcite{\IncredibleFormula}:

  \state Corollary
  \label \ReformulationOfIncredibleFormula
  If $a,b\in\Na$ then
  $a^*a\rc b^*b$ if and only if
  $\zeta(a)\zeta(b)^*\in\fullcp.$

  We do not claim that the relation of being relatively continuous is
reflexive, symmetric, or transitive.  Instead we have:

  \state Proposition
  \label \RCProperties
  Let $a$, $b$, and $c$ be $\a$-integrable elements.
  \izitem If $a\rc c$ and $b\rc c$ then $a+b\rc c$.
  \zitem If $a\rc b$ then $b^*\rc a^*$.
  \zitem If $a\rc b$ then for every $w \in \Gdual$ and every
$m\in\Mult_w(A)$ one has that both
  $ma\rc b$ and $a\rc bm$.
  \zitem If $a\rc c$ and $c\rc b$ then
  for every $w \in \Gdual$ we have that
  $a\F cw\rc b$ and $a\rc \F cwb$.

  \proof The first assertion is trivial.  In order to prove (ii) note
that
  by \lcite{\FormulasForFourier}.(i) we have for all $x,y,z\in\Gdual$
that
  $$
  \big(\F{b^*}{xz} \F{a^*}y - \F{b^*}x \F{a^*}{zy}\big)^* =
  - \big(
  \F {a}{y\inv z\inv} \F {b}{x\inv} -
  \F {a}{y\inv} \F {b}{z\inv x\inv}
  \big),
  $$
  from where (ii) follows easily.
  As for (iii), using \lcite{\FormulasForFourier}.(ii), we have for
all $x,y,z\in\Gdual$ that
  $$
  \[\F {ma}{xz} \F{b}{y} - \F{ma}{x} \F{b}{zy}\] =
  \[m\F {a}{w\inv xz} \F{b}{y} - m\F{a}{w\inv x} \F{b}{zy}\] \$\leq
  \[m\] \[\F {a}{w\inv xz} \F{b}{y} - \F{a}{w\inv x} \F{b}{zy}\].
  $$
  This proves that $ma\rc b$.  The proof that $a\rc bm$ goes along
similar lines.
  With respect to (iv) observe that by
\lcite{\FormulasForFourier}.(iv) we have
  $$
  \[\F {a\F cw}{xz} \F{b}{y} - \F{a\F cw}{x} \F{b}{zy}\] \$=
  \[\F {a}{xzw\inv} \F {c}{w}\F{b}{y} - \F{a}{xw\inv}\F cw \F{b}{zy}\]
  \$\leq
  \[\F {a}{xzw\inv} \F {c}{w}\F{b}{y} - \F {a}{xw\inv} \F
{c}{zw}\F{b}{y}\]
  \$+ \[\F {a}{xw\inv} \F {c}{zw}\F{b}{y} - \F{a}{xw\inv}\F cw
\F{b}{zy}\]
  \$\leq
  \[\F {a}{xzw\inv} \F {c}{w} - \F {a}{xw\inv} \F {c}{zw}\] \[b\]_1
  + \[a\]_1\[\F {c}{zw}\F{b}{y} - \F cw \F{b}{zy}\],
  $$
  where we have used \lcite{\NormOfFourierVsOneNorm} in the last step.
This proves that $a\F cw\rc b$.  Again the proof that $a\rc \F cwb$
follows similarly.
  \proofend

  Relative continuity enjoys a curious hereditary property which we
discuss next.

  \begingroup
  \def\aa{a}
  \def\b{b}
  \def\c{c}
  \def\d{a_1}
  \def\e{b_1}
  \def\f{c_1}
  \state Proposition
  \label \Herege
  Let $\aa,\b,\c\in A$ be positive $\a$-integrable elements such that
  $\aa\leq \b\rc \c$.
  Then $\aa\rc \c$.

  \proof
  Write $\aa=\d^*\d$, $\b=\e^*\e$, and $\c=\f^*\f$, for $\d,\e,\f\in
A$.  Since $\aa\leq \b$ we have for all $v$ in $\H$ that
  $\[\d(v)\] =
  \<\d^*\d(v),v\>^{1/2}\leq \<\e^*\e(v),v\>^{1/2} =
  \[\e(v)\]$.  It therefore follows that there exists a bounded
operator $T$ on $\H$ with $\[T\]\leq1$ such that $\d=T\e$.

Even though $T$ may not be in $A$ we would like to make sense of the
expression $\pi(T)$ and also to show that $\pi(T)\zeta(\e) =
\zeta(T\e)$, mimicking \lcite{\NewActionOfCrossedProduct}.(i).  In
order to accomplish this let us suppose, without loss of generality,
that there is a strongly continuous unitary representation $U$ of $G$
on $\H$ such that $\a_t(a)=U_t a U_t\inv$, for all $a$ in $A$.
  We may then define
  $$
  \pi : \B(\H) \to \B(L_2(G,\H)),
  $$
  by
  $$
  \pi(S)\xi\calcat t = U_t\inv S U_t \xi(t)
  \for S\in\B(\H)\for \xi\in L_2(G,\H)\for t\in G,
  $$
  extending the representation $\pi$ of section
\lcite{\BasicFactsSection}.  Observe that $\pi(S)\xi$ is a measurable
function on $G$, and thus represents an element of $L_2(G,\H)$,
because $U$ is strongly continuous.
  We now claim that
  $\pi(T)\zeta(\e) = \zeta(\d)$.  In fact, given $v\in\H$ and $t\in G$
we have
  $$
  \zeta(\d)v\calcat t =
  \a_t\inv(T\e)v =
  U_t\inv T\e U_t v =
  U_t\inv T U_t \a_t\inv(\e) v \$=
  U_t\inv T U_t \(\zeta(\e) v\calcat t\) =
  \Big(\pi(T)\zeta(\e) v\Big)\calcat t,
  $$
  thus proving our claim.  In order to show that $\aa\rc \c$
  let $x,y,z\in \Gdual$ and observe that by \lcite{\FourierViaV}
  $$
  \[\F \aa{xz}\F \c y - \F \aa x\F \c {zy}\] \$=
  \[\zeta(\d)^*V_{xz}\zeta(\d) \zeta(\f)^*V_{y}\zeta(\f) -
    \zeta(\d)^*V_{x}\zeta(\d) \zeta(\f)^*V_{zy}\zeta(\f)\] \$\leq
  \[\zeta(\d)^*V_x\]
  \[V_z\pi(T)\zeta(\e) \zeta(\f)^*-
    \pi(T)\zeta(\e) \zeta(\f)^*V_z\]
    \[V_{y}\zeta(\f)\] = \ldots
  $$
  It is easy to show that $\pi(T)$ commutes with $V_z$.  Therefore the
above equals
  $$
  \ldots =
  \[\zeta(\d)\]
  \[\pi(T)\zeta(\e) \zeta(\f)^*-
    \pi(T)V_z\inv\zeta(\e) \zeta(\f)^*V_z\] \[\zeta(\f)\] \$\leq
  \[\zeta(\d)\]
  \[\pi(T)\]
  \[\zeta(\e) \zeta(\f)^*-V_z\inv\zeta(\e) \zeta(\f)^*V_z\]
\[\zeta(\f)\].
  $$
  Since, by hypothesis we have that $\e^*\e\rc\f^*\f$ we conclude by
\lcite{\ReformulationOfIncredibleFormula} that
$\zeta(\e)\zeta(\f)^*\in\fullcp$ and hence that this is a
$V$-continuous operator.  It follows that the last expression
displayed above tends to zero as $z\to e$ and hence that
  $\aa\rc c$.
  \proofend
  \endgroup

  \state Proposition
  \label \IntIsRelContinuous
  Let $a$ and $b$ be $\a$-integrable elements such that $a\rc b$.
Then, for each $g\in C_c(G)$ one has that $a'\rc b$, where $a'=\int_G
g(t) \a_t(a) \d t$.

  \proof
  Given $x,y,z\in\Gdual$ we have, using
\lcite{\FourierTransfOfConvolution}, that
  $$
  \[\F {a'}{xz} \F by - \F {a'}x \F b{zy}\] =
  \[\^g\(xz\) \F a{xz} \F by - \^g\(x\)\F ax \F b{zy}\] \$\leq
  \[\^g\(xz\) \F a{xz} \F by - \^g\(x\) \F a{xz} \F by\] + \[\^g\(x\)
\F a{xz} \F by - \^g\(x\)\F ax \F b{zy}\] \$\leq
  \left|\^g\(xz\) - \^g\(x\)\right| \[a\]_1 \[b\]_1 +
\left|\^g\(x\)\right|\[ \F a{xz} \F by - \F ax \F b{zy}\],
  $$
  which tends to zero, uniformly in $x$ and $y$, as $z\to e$, because
$\^g$ is uniformly continuous and bounded.
  \proofend

  We will often be concerned with sets of mutually relatively
continuous elements.  For this reason we make the following:

  \definition
  \label \RelativelyContinuousSpectrallyInvariant
  Let $\W$ be a set of $\a$-integrable elements.
  \izitem We say that $\W$ is a \stress{relatively continuous} set if
for every $a,b\in \W$ one has that $a\rc b$.
  \zitem We say that $\W$ is \stress{spectrally invariant} if, given
$a\in \W$ and $x\in\Gdual$, one has that $\F ax \W\subseteq \W$ and
  $\W\F ax \subseteq \W$.
  \zitem The \stress{spectrally invariant hull} of $\W$, denoted
$\widetilde \W$, is the intersection of all spectrally invariant sets
of $\a$-integrable elements containing $\W$.

  Observe that the intersection of any number of spectrally invariant
sets is again spectrally invariant and so is $\widetilde \W$.

We should also note that the elements of a relatively continuous set
satisfy
  $a\rc a$, which is by no means automatic.

  \state Proposition
  \label \ConsertaMasNaoEstraga
  Let $\W$ be a relatively continuous set of $\a$-integrable elements.
Then the spectrally invariant hull $\widetilde \W$ of $\W$ is also
relatively continuous.

  \proof
  In order to avoid repetition, during the course of this proof $a$,
$b$, and $c$, with or without subscripts, will always refer to
elements of $\W$, and $x$, $y$, and $z$ will denote elements of
$\Gdual$.

Given $n$-vectors
  $\vec a = (a_1,\ldots,a_n)$ and $\vec x = (x_1,\ldots,x_n)$ let
  $$
  \F{\vec a}{\vec x} = \F {a_1}{x_1}\ldots\F {a_n}{x_n}.
  $$ If $\vec b = (b_1,\ldots,b_m)$ and $\vec y = (y_1,\ldots,y_m)$ is
another such pair of vectors consider the set
  $\F{\vec a}{\vec x} \W \F{\vec b}{\vec y}$.
  Its elements are thus of the form
  $u=\F{\vec a}{\vec x} c \F{\vec b}{\vec y}$.  Observe that, by
\lcite{\FormulasForFourier}, $u$ is $\a$-integrable and for
$z\in\Gdual$ we have
  $$
  \F uz =
  \F{\vec a}{\vec x}
    \F c{z'}
  \F{\vec b}{\vec y},
  $$
  where $z'=(x_1\ldots x_n)\inv z (y_1\ldots y_m)\inv$.
  It is thus clear that the union of all sets of the form $\F{\vec
a}{\vec x} \W \F{\vec b}{\vec y}$, as above, is spectrally invariant
and hence, being the smallest among those which contain $\W$,
coincides with $\widetilde \W$.

It therefore suffices to show that any two elements
  $$
  u=\F{\vec a}{\vec x} c \F{\vec b}{\vec y}
  \and
  u'=\F{\vec {a'}}{\vec {x'}} c' \F{\vec {b'}}{\vec {y'}}
  $$
  of the above form satisfy $u\rc u'$.  In view of
\lcite{\RCProperties}.(iii) it suffices to consider the case where
  $$
  u= c \F{\vec b}{\vec y}
  \and
  u'=\F{\vec {a'}}{\vec {x'}} c'.
  $$
  In order to prove this we use induction on $m=|\vec b|+|\vec {a'}|$,
where $|{\cdot}|$ denotes the number of coordinates of a given vector.
If $m=0$ this follows from the hypothesis.  Now, suppose that $m>0$
and hence either $|\vec b|>0$ or $|\vec a'|>0$.  Suppose first that
$n:= |\vec b|>0$ and let $\vec {b^\sharp} = (b_1,\ldots,b_{n-1})$ and
  $\vec {y^\sharp} = (y_1,\ldots,y_{n-1})$ so that
  $\F {\vec b}{\vec y} = \F {\vec {b^\sharp}}{\vec {y^\sharp}} \F
{b_n}{y_n}$.  By the induction hypothesis we have that
  $$
  c\F{\vec {b^\sharp}}{\vec {y^\sharp}} \rc
  b_n \rc
  u'.
  $$
  By \lcite{\RCProperties}.(iv) we then conclude that $u\rc u'$.  If
$|\vec b|=0$ then $|\vec a'|>0$ and a similar argument may be used to
complete the proof.
  \proofend

  \definition
  \label \AlgSet
  For any subset $X\subseteq A$ we denote by $\Alg X$ the smallest
$\a$-invariant sub-\cstar-algebra of $A$ containing $X$.

  \state Proposition
  \label \SIHHasSameAlgebra
  If\/ $\W\subseteq A$ is a set of $\a$-integrable elements then
$\Alg\W=\Alg{\widetilde \W}$, where
  $\widetilde \W$ is the spectrally invariant hull of $\W$.

  \proof
  Let $B$ be an $\a$-invariant sub-\cstar-algebra of $A$.  We claim
that the set $B_1$ of $\a$-integrable elements in $B$ is spectrally
invariant.  In fact, given $a,b\in B_1$, $x\in\Gdual$, and $t\in G$ we
have that $\du xt \a_t(a)b\in B$ and hence the unconditional integral
of this expression with respect to $t$, namely $\F axb$, also belongs
to $B$.  Since $\F axb$ is $\a$-integrable by
\lcite{\FormulasForFourier}.(ii) we have that $\F axb\in B_1$.
Similarly one shows that $b\F ax\in B_1$, hence proving our claim.

  Taking $B=\Alg\W$ in the argument above we thus conclude that
$\widetilde \W\subseteq \Alg\W$, from which it follows that
$\Alg{\widetilde \W}\subseteq \Alg\W$.  The converse inclusion is
trivial.
  \proofend

By definition a \stress{cone} is a subset of $A$ which is closed under
addition and multiplication by positive (i.e $\geq 0$) scalars.
  We would now like to study relative continuity for cones.  Recall
that $\Pa$ denotes the hereditary cone of all positive $\a$-integrable
elements of $A$.

  \state Proposition
  \label \PropertiesOfMaximalCones
  Let $\P\subseteq\Pa$ be a relatively continuous cone which is
maximal among the relatively continuous cones contained in $\Pa$.
Then
  \izitem $\P$ is hereditary.
  \zitem If $g\in C_c(G)$ is a positive function and $a\in\P$ then
  $a'\in\P$, where
  $a' = \int_G g(t)\a_t(a) \d t$.

  \proof
  It is useful to keep in mind that relative continuity is a symmetric
relation among self-adjoint elements by \lcite{\RCProperties}.(ii).

  Suppose that $0\leq a\leq b$, where $a\in A$ and $b\in \P$.  By
\lcite{\BasicPropertiesOfModule}.(i) we have that $a$ is
$\a$-integrable.  Define
  $\widetilde\P=\{\lambda a + p : \lambda \geq 0, \ p\in\P\}$, which
then turns out to be a cone consisting of $\a$-integrable positive
elements.  We claim that $\widetilde\P$ is relatively continuous.  In
order to prove this it clearly suffices to show that
  $a\rc a$ and that $a\rc p$ for all $p\in\P$.
  Since $\P$ is relatively continuous we have that $a\leq b\rc p$ for
every $p\in\P$.  So $a\rc p$ by \lcite{\Herege}.  It now follows that
$a\leq b\rc a$ and hence that $a\rc a$.  This proves that
$\widetilde\P$ is relatively continuous.  Since $\P$ is maximal by
hypothesis, we conclude that $\widetilde\P=\P$ and hence that
$a\in\P$.

To prove the second assertion let
  $\widetilde\P=\{\lambda a' + p : \lambda \geq 0, \ p\in\P\}$, which
is a subcone of $\Pa$ by \lcite{\IntIsIntegrable}.  Again we claim
that $\widetilde\P$ is relatively continuous and to prove it we must
once more show that
  $a'\rc a'$ and that $a'\rc p$ for all $p\in\P$.  Since $a\in \P$ we
have that
  $a\rc p$ and hence
  \lcite{\IntIsRelContinuous} implies that $a'\rc p$.  In particular,
taking $p=a$ we have that $a\rc a'$ and
  \lcite{\IntIsRelContinuous} finally gives $a'\rc a'$.  This proves
that $\widetilde\P$ is relatively continuous and the conclusion
follows as above.
  \proofend

  \state Corollary
  \label \TheseAnimalsExist
  Let $\W$ be a relatively continuous set of positive $\a$-integrable
elements.  Then $\W$ is contained in a relatively continuous cone $\P$
satisfying \lcite{\PropertiesOfMaximalCones}.(i-ii).

  \proof
  It is easy to see that the set $\P_0$ formed by all linear
combinations with positive coefficients of elements of $\W$ is a
relatively continuous cone.  By Zorn's Lemma let $\P$ be a maximal
relatively continuous cone containing $\P_0$.  Then $\P$ satisfies the
required properties by \lcite{\PropertiesOfMaximalCones}.
  \proofend

  \section {{\MR} equivalence}
  In this section we will apply the knowledge gained so far to the
problem of establishing a {\MR} equivalence between subalgebras of
$\Mult_e(A)$ and ideals of $\fullcp$.

  For the time being we shall fix a relatively continuous cone
$\P\subseteq\Pa$ satisfying \lcite{\PropertiesOfMaximalCones}.(i-ii).
Let
  $$
  \N=\{a\in A : a^*a\in \P\}
  \and
  \M =\span (\P).
  $$
  Since $\P$ is supposed hereditary we have by
\lcite{\ConeModuleAlgebra} that $\N$ is a left ideal of $A$, clearly
contained in $\Na$, and $\M$ is a subalgebra of $A$ contained in
$\Ma$.

We shall denote by $\X$ the closure of $\zeta(\N)$ in $\Xa$.

  \state Proposition
  One has that
  $\(\fullcp\)\X \subseteq \X$.

  \proof
  We will obviously use \lcite{\LeftCrossedProductModule}, observing
that the second hypothesis required there holds for ${\cal D}=A$ since
$\N$ is a left ideal.  As for \lcite{\LeftCrossedProductModule}.(i)
let $a'=\int_G g(t) \a_t(a) \d t$ where $g\in C_c(G)$ and $a\in \N$.
Observe that, by \lcite{\Ineq}, we have
  $$
  a'^*a' \leq
  |\text{supp}(g)| \int_G |g(t)|^2 \a_t(a^*a) \d t.
  $$
  Since $|g(t)|^2$ is a positive compactly supported function on $G$
and $a^*a\in\P$ we have that the right hand side above gives an
element in $\P$ (by the fact that $\P$ satisfies
\lcite{\PropertiesOfMaximalCones}.(ii)).  Also, since $\P$ is
hereditary we have that $a'^*a'\in\P$ and hence that $a'\in\N$.  This
shows that $\N$ satisfies the requirements in
\lcite{\LeftCrossedProductModule} and hence the proof is complete.
  \proofend

  \state Theorem
  \label \AbstractMR
  Let $\P$ be a relatively continuous cone consisting of positive
$\a$-integrable elements of $A$.  Suppose that $\P$ satisfies
\lcite{\PropertiesOfMaximalCones}.(i-ii) (e.g. if $\P$ is maximal).
Set
  $\N=\{a\in A : a^*a\in \P\}$ and $\M =\span (\P)$, and let $\X$ be
the closure of $\zeta(\N)$ within $\Xa$.  Then
  $\X$ is a left Hilbert $\(\fullcp\)$--module and hence it
establishes a {\MR} equivalence between $\overline{\X^*\X}$, which is
a subalgebra of $\Mult_e(A)$ coinciding with $\overline{\F {\M}e}$,
and the ideal $\overline{\X\X^*}$ in $\fullcp$.

  \proof
  By the last Proposition we see that $\X$ is a left
$\(\fullcp\)$--module under multiplication of operators.
  Given $T,S\in\X$ we claim that $TS^*\in\fullcp$.  In fact, by
definition of $\X$ it is enough to show that $\zeta(a)\zeta(b)^*\in
\fullcp$ for all $a,b\in\N$.
  Since $a^*a$ and $b^*b$ belong to $\P$, and $\P$ is relatively
continuous, we have that $a^*a\rc b^*b$ and hence the claim follows
{}from \lcite{\ReformulationOfIncredibleFormula}.  Therefore $\X$ is a
left Hilbert $\fullcp$--module for the inner-product $(T|S)_L :=
TS^*$, and hence $\overline{\X\X^*}$ is an ideal in $\fullcp$.

By \lcite{\ColumnOperators}.(iii) we have that
  $\zeta(a)^*\zeta(b) = \F {a^*b}e$ and hence, using
\lcite{\ConeModuleAlgebra}.(ii), we have
  $$
  \zeta(\N)^* \zeta(\N) = \F {\N^*\N} e = \F {\M} e \subseteq
\Mult_e(A).
  $$
  It follows that
  $\overline{\X^*\X} =
  \overline{\F{\M}e}.$
  \proofend

  Using our last result in combination with \lcite{\TheseAnimalsExist}
we obtain:

  \state Corollary
  Let $\W$ be a relatively continuous set of positive $\a$-integrable
elements.  Then there exists a {\MR} equivalence between a subalgebra
of $\Mult_e(A)$ containing $\F {\W}e$ and an ideal in $\fullcp$.

  This brings about the following:

  \sysstate{Question}{\rm}
  {\label \DensityOfRCCone Suppose that $\Pa$ is dense in $A_+$.  Does
if follow that there exists a relatively continuous cone $\P$ which is
also dense in $A_+$?}

  \sysstate{Question}{\rm}
  {\label \QuestionOnUFA Is the subalgebra $\overline{\X^*\X}$ of
$\Mult_e(A)$ always the same for different maximal relatively
continuous cones?}

If the above question is answered affirmatively then one should
probably call $\overline{\X^*\X}$ the ``{\GFPA}'' of $\a$.  We shall
have more to say about this question in what follows.

  \section {Dual action on Fell bundles}
  \label \SectionOnDualAction
  We would now like to discuss an important example of
\cstar-dynamical system, namely that of the dual action on the
cross-sectional \cstar-algebra of a Fell bundle.

This action is particularly well behaved from the point of view we
wish to discuss and, in particular, question \lcite{\DensityOfRCCone}
will be answered affirmatively.

Another specially important aspect of this example is that the concept
of spectral subspace is built in, as the bundle fibers.  For this
reason we will take this example as the model for our spectral theory
to be developed in the next section.

Let $G$ be an abelian group with Pontryagin dual $\Gdual$ and let
$\B=\{B_x\}_{x\in\Gdual}$ be a Fell bundle over $\Gdual$ (see
\cite{\FD} for a comprehensive treatment of the theory of Fell
bundles, also referred to as \cstar-algebraic bundles).  Let $C^*(\B)$
be the cross sectional \cstar-algebra of $\B$ \scite{\FD}{VIII.17.2}
defined to be the enveloping \cstar-algebra of the Banach *-algebra
$L_1(\B)$ formed by the integrable sections \scite{\FD}{VIII.5.2}.

As before we will fix a faithful representation of $C^*(\B)$ on a
Hilbert space $\H$ and hence will think of $C^*(\B)$ as an algebra of
operators on $\H$.

Denote by $C_c(\B)$ the dense sub-algebra of $L_1(\B)$ formed by the
continuous, compactly supported sections \scite{\FD}{II.14.2}. We
remark that our notation differs from \cite{\FD} with respect to
$C_c(\B)$.

As described in \scite{\FD}{VIII.5.8} (see also
\scite{\unconditional}{Section 5}) each $b_x$ in each fiber $B_x$ of
$\B$ defines a multiplier of $C^*(\B)$ such that
  $$
  (b_x f) (y) = b_x f(x\inv y)
  \and
  (f b_x) (y) = f(y x\inv) b_x,
  $$
  for all $f\in L_1(\B)$ and $y\in \Gdual$.

The \stress{dual action} of $G$ on $C^*(\B)$, which we denote simply
by $\a$, is determined by the expression
  $$
  \a_t(f)\calcat x = \idu xt f(x),
  $$
  for all $f\in L_1(\B)$, $t\in G$, and $x\in\Gdual$.  See
\scite{\unconditional}{Section 5} for details but beware of the change
in conventions which we attempted to explain in footnote
\lcite{\Convention}.
  It is easy to see that the natural extension of $\a$ to the
multiplier algebra $\Mult(C^*(\B))$ satisfies
  $$
  \a_t(b_x) = \idu xt b_x,
  $$
  for all $b_x$ in any $B_x$.  Therefore $B_x\subseteq
\Mult_x(C^*(\B))$ for all $x$ in $\Gdual$.  Since this inclusion may
be proper, one may argue over which set deserves to be called the
``$x$-spectral subspace'' for the dynamical system $(C^*(\B),G,\a)$.
We believe that $B_x$ is the ``correct'' choice and we hope to
convince the reader of this in what follows.  In particular we propose
defining the ``{\GFPA}'' to be $B_e$ among the many subalgebras of
$\Mult_e(C^*(\B))$ which occur in examples.

An argument in favor of this point of view is related to the Takai
duality Theorem:
  given an action of $\Gdual$ on a \cstar-algebra $A$ consider the
corresponding semi-direct product bundle $\B$ \scite{\FD}{VIII.4.2}.
Its cross sectional algebra is $A\crossproduct \Gdual$ and the dual
action, as defined above, coincides with the classical notion of dual
action \scite{\Ped}{7.8.3}, up to a sign convention.  Moreover the
crossed product by the dual action is isomorphic to
$A\*\K(L_2(\Gdual))$ \scite{\Ped}{7.9.3}.  Thus, if we want the
``{\GFPA}'' to be {\MR} equivalent to the crossed product (see
\cite{\RieffelProper} and \scite{\RieffelIntegrable}{Section 6}) a
natural choice is to take it to be $A$, which is precisely the unit
fiber of the semi-direct product bundle!

According to Theorem 5.5 of \cite{\unconditional} every element of
$C^*(\B)$ of the form $p=f^* f$, with $f\in C_c(\B)$, is
$\a$-integrable and in
  addition\footnote{\footcntr}
  {\eightpoint Observe that \cite{\unconditional} uses a different
sign convention both for the dual action and for the Fourier
coefficients.  Nevertheless these differences compensate each other in
such a way that the statement of Theorem 5.5 in \cite{\unconditional}
remains valid as it stands.},
  $$
  \F px = p(x),
  \leqno{(\seqnumbering)}
  \label \FourierTransformInBundle
  $$
  where the term in the right hand side is to be interpreted as a
multiplier of $C^*(\B)$ as explained above.  It then follows that
$C_c(\B)$ is contained in $\Na$.

  \state Proposition
  \label \RelContInBundles
  The subset $C_c(\B)^2\subseteq C^*(\B)$ of all linear combinations
of products of elements in $C_c(\B)$ is a relatively continuous set of
$\a$-integrable elements.

  \proof
  As seen above $a^*a$ is $\a$-integrable for all $a\in C_c(\B)$.  By
the polarization formula it is then easy to see that $a^*b$ is
$\a$-integrable for all $a,b\in C_c(\B)$.  Since $C_c(\B)$ is
self-adjoint we then conclude that any element of $C_c(\B)^2$ is
$\a$-integrable.

  We now have to prove that $a^*b\rc c^*d$ for all $a,b,c,d\in
C_c(\B)$.  Again by the polarization formula it is enough to verify
the case $a^*a\rc b^*b$.
  In order to do this we shall use the implication
(iii)$\Rightarrow$(ii) of \lcite{\IncredibleFormula}.  So let $p=a^*
a$ and $q=b^* b$ and observe that by \lcite{\FourierTransformInBundle}
we have
  $$
  \lim_{x\to e}\[\F px \E q - \E p \F qx\] =
  \lim_{x\to e}\[p(x) q(e) - p(e) q(x)\],
  $$
  and
  $$
  \lim_{x\to e}\[\F px \F q{x\inv} - \E p \E q\] =
  \lim_{x\to e}\[p(x)q(x\inv) - p(e)q(e)\].
  $$
  Since the norm and multiplication are continuous on $\B$, and both
$p$ and $q$ are continuous sections, we see that the limits above
vanish.  This shows that \lcite{\IncredibleFormula}.(iii) holds and
hence also \lcite{\IncredibleFormula}.(ii), concluding the proof.
  \proofend

As a consequence we have:

  \state Corollary
  \label \LipInCompactSupport
  Let $a,b\in C_c(\B)$.  Then $\zeta(a)\zeta(b)^*$ belongs to
$C^*(\B)\crossproduct_\a G$.

  \proof
  By the previous result we have that $a^*a\rc b^*b$.  So the
conclusion follows from \lcite{\ReformulationOfIncredibleFormula}.
  \proofend

  Observe that $C_c(\B)^2\cap C^*(\B)_+$ is dense in $C^*(\B)_+$ and
hence $C^*(\B)_+$ contains a dense cone which is relatively
continuous.  This answers question \lcite{\DensityOfRCCone}
affirmatively.

There is a very natural Hilbert module associated to any given Fell
bundle.  That is the $B_e$-module, denoted $L_2(\B)$, defined to be
the completion of $C_c(B)$ under the $B_e$--valued inner-product
defined by the integral $\int_{\Gdual} f(x)^*g(x)\d x$, for $f$ and
$g$ in $C_c(\B)$.  Our next result is intended to relate it to the
present situation.  Recall that $\Rip fg = \sint_G \a_t(f^*g) \d t$,
as defined in \lcite{\FirstOcurenceOfRip}.

  \state Proposition
  \label \OldAndNewInnerProd
  Let $f,g\in C_c(B)$.  Then
  $\Rip fg = \int_\Gdual f(x)^*g(x) \d x$.

  \proof
  By the polarization identity it suffices to consider $f=g$.  In this
case we have
  $$
  \Rip ff =
  \sint_G \a_t(f^* f) \d t =
  \F {f^* f}e \={\FourierTransformInBundle}
  (f^* f)(e) =
  \int_\Gdual f(x)^*f(x) \d x.
  \proofend
  $$

  We therefore conclude that the inclusion $C_c(\B)\subseteq \Na$
respects the corresponding $B_e$ valued inner-products and hence that
$L_2(\B)$ is isomorphic to the closure of $\zeta(C_c(\B))$ within
$\Xa$, as Hilbert $B_e$--modules.

  \state
  \label \CompactSupportisModule
  Proposition
  Let $\X$ be the subset of $\Xa$ obtained by closing
$\zeta(C_c(\B))$.  Then
  $$
  \big(C^*(\B)\crossproduct_\a G\big) \X \subseteq \X.
  $$

  \proof
  We will derive this from \lcite{\LeftCrossedProductModule},
observing that \lcite{\LeftCrossedProductModule}.(ii) is satisfied for
${\cal D} = C_c(\B)$.  In order to prove
\lcite{\LeftCrossedProductModule}.(i),
  let $a\in C_c(\B)$ and denote the closed support of $a$ by $K$.
Denote by $C_K(\B)$ the subset of $C_c(\B)$ formed by the continuous
sections of $\B$ vanishing on $G\setminus K$.  Equipped with the
supremum norm $C_K(\B)$ is a normed space and in fact a Banach space
by \scite{\FD}{II.13.13}.
  Consider the map
  $$
  \psi: t\in G \mapsto \a_t(a) \in C_K(\B),
  $$
  (that $\a_t(a)$ is in $C_K(\B)$ follows from \scite{\FD}{II.13.14}).
We claim that $\psi$ is continuous.  In fact, if a net $\{t_i\}$
converges to $t$ in $G$, then
  $$
  \lim_i \sup_{x\in K} \[\tdu xt - \tdu x{t_i}\] = 0
  $$
  by the Pontryagin--van Kampen duality Theorem \scite{\HR}{24.8},
{}from which the claim follows easily.  Given $g$ in $C_c(G)$ we then
have that the Bochner integral
  $$
  \int_G g(t) \a_t(a) \d t
  $$
  converges and hence defines an element $a'\in C_K(\B)$.  Since the
inclusions
  $$
  C_K(\B) \to L_1(\B) \to C^*(\B)
  $$
  are obviously continuous we conclude that the above integral, if
seen as the integral of a $C^*(\B)$--valued function, also converges
to $a'$.
  It therefore follows that $C_c(\B)$ satisfies the hypothesis of
\lcite{\LeftCrossedProductModule}, concluding the proof.
  \proofend

  The following result is a generalization of the Takai duality for
the kind of actions we are dealing with.  It is perhaps also an
argument in favor of defining the {\GFPA} to be the unit fiber
algebra, in the case of a dual action.

  \state Theorem
  \label \MoritaForFellBundles
  Let $\X$ be the closure of $\zeta(C_c(\B))$ as above.  Then
  $\X$ is a Hilbert bimodule over the unit fiber algebra $B_e$ on the
right hand side and $C^*(\B)\crossproduct_\a G$ on the left.  Moreover
$\overline{\X\X^*}$ is an ideal in $C^*(\B)\crossproduct_\a G$ which
is {\MR} equivalent to $B_e$.  The imprimitivity bimodule implementing
this equivalence may be taken to be $\X$.

  \proof
  As seen above $C_c(\B)\subseteq \Na$ and $B_e\subseteq
\Mult_e(C^*(\B))$.  Using \lcite{\ColumnOperators}.(ii)
  we then have that
  $\zeta\big(C_c(\B)\big)B_e \subseteq \zeta\big(C_c(\B)\big)$ and
hence that
  $\X B_e \subseteq \X$.  In other words $\X$ is a right
$B_e$--module.  From \lcite{\CompactSupportisModule} we have that
$\big(C^*(\B)\crossproduct_\a G\big) \X \subseteq \X$ and hence $\X$
is a left $\big(C^*(\B)\crossproduct_\a G\big)$--module.

  For $T,S\in\X$ consider the inner-products
  $$
  (T|S)_R = T^*S
  \and
  (T|S)_L = TS^*.
  $$
  By \lcite{\ColumnOperators}.(iii) and \lcite{\OldAndNewInnerProd} we
have that the range of $(\cdot|\cdot)_R$ is contained in $B_e$, and by
\lcite{\LipInCompactSupport}, that the range of $(\cdot|\cdot)_L$ is
contained in $C^*(\B)\crossproduct_\a G$.  Thus $\X$ is a Hilbert
$\big(C^*(\B)\crossproduct_\a G\big)$--$B_e$--bimodule.

We claim that $\overline{\X^*\X}$ coincides with $B_e$.  In fact,
given a positive element $b\in B_e$ choose by \scite{\FD}{Appendix C}
a continuous section $f$ of $\B$ such that $f(e)=b^{1/2}$.  Given
$\varepsilon>0$ let $\Omega$ be a neighborhood of $e$ in $\Gdual$ such
that $x\in\Omega$ implies that
  $\[f^*(x)f(x)-b\]<\varepsilon$.  Also take $g\in C_c(\Gdual)$ such
that $\text{supp}(g)\subseteq\Omega$ and $\int_\Gdual |g(x)|^2 \d x
=1$.  It follows that the section $f'$, given by the pointwise product
$f'=gf$, is in $C_c(\B)$ and that
  $$
  \[\int_\Gdual f'(x)^*f'(x) \d x - b\]<\varepsilon.
  $$
  On the other hand we have that
  $$
  \int_\Gdual f'(x)^*f'(x) \d x \={\lcite{\OldAndNewInnerProd}}
  \Rip{f'}{f'} \={\ColumnOperators}
  \zeta(f')^*\zeta(f')\in \X^*\X,
  $$
  proving that $b$ is in $\overline{\X^*\X}$.  This proves our claim.
The remaining statements are routine.
  \proofend

  \section {Spectral theory}
  As seen above, dual actions in the context of Fell bundles gives
rise to well behaved \cstar-dynamical systems from the point of view
of the questions we propose to treat.
  It is therefore desirable to characterize the abelian group actions
which arise as such.  Another important reason why one would like to
describe an action as a dual action is the classification Theorem for
Fell bundles \scite{\TPA}{7.3} according to which every Fell bundle is
stably isomorphic to the semi-direct product bundle for a twisted
partial action of the base group on the unit fiber algebra.

  Precisely speaking, given a strongly continuous action $\a$ of an
abelian group $G$ on a \cstar-algebra $A$, satisfying suitable
hypothesis, we wish to find a Fell bundle $\B$ over the dual group
$\Gdual$ such that $A$ is isomorphic to $C^*(\B)$ under an isomorphism
which puts $\a$ in correspondence with the dual action of $G$ on
$C^*(\B)$.  Recall from \lcite{\RelContInBundles} that dual actions
contain a dense relatively continuous set of integrable elements.  It
is therefore natural to require the existence of relatively continuous
sets, which we do next.  For convenience we have chosen to require
certain conditions, namely (iv) and (v) below, which will eventually
be removed when we state the last and main Theorem of this section.

  \sysstate{Hypothesis}{\rm}
  {\label \WorkingHypothesis
  We shall assume, until further notice, that we are given a subset
$\W\subseteq A$ such that
  \izitem $\W^*=\W$,
  \zitem $\W$ consists of $\a$-integrable elements,
  \zitem $\W$ is relatively continuous,
  \zitem $\W$ is spectrally invariant, and
  \zitem $\W$ is a linear subspace of $A$.
  }

   Our main goal will be to construct a Fell bundle $\B$ over $\Gdual$
whose cross-sectional algebra $C^*(\B)$ is isomorphic to $\Alg\W$
(Definition \lcite{\AlgSet}) under an isomorphism which is covariant
for the dual action of $G$ on $C^*(\B)$ and the action $\a$ on $A$.

  \definition
  For each $x$ in $\Gdual$ let $B_x$ be defined to be the closure of
the set
  $$
  \{\F {a}x \colon a\in\W\}
  $$
  within $\Mult_x(A)$.

  Since $\W$ is a linear space, it is clear that $B_x$ is a Banach
space.

  \state Proposition
  For any $x,y\in\Gdual$ one has that
  $B_x B_y\subseteq B_{xy}$ and $B_x^* = B_{x\inv}$.

  \proof
  Given $a,b\in\W$ we have by \lcite{\FormulasForFourier}.(iv) that
  $$
  \F {a}x \F {b}y =
  \F {a \F {b}y}{xy}.
  $$
  Since $\W$ is spectrally invariant we have that $a \F {b}y$ is in
$\W$ and hence that
  $\F {a}x \F {b}y$ belongs to $B_{xy}$.  This proves the first
statement.  The second statement follows immediately from by
\lcite{\FormulasForFourier}.(i) and the assumption that $\W$ is
self-adjoint.
  \proofend

  So we see that the collection $\B = \{B_x\}_{x\in\Gdual}$ forms a
Fell bundle over the group obtained by replacing the topology of
$\Gdual$ with the discrete topology.  In order to make $\B$ into a
Fell bundle over $\Gdual$ (with its own topology) we will use
\scite{\FD}{II.13.18}.  That is, we must provide a linear space
$\Gamma$ of sections of $\B$ such that
  \iItem{(a)} for each $f\in\Gamma$ the numerical function $x\mapsto
\[f(x)\]$ is continuous on $\Gdual$, and
  \Item{(b)} for each $x$ in $\Gdual$ the set $\{f(x)\colon
f\in\Gamma\}$ is dense in $B_x$.

  \state Proposition
  \label \TheBanachBundle
  Let $\Gamma$ be the linear space of sections $f$ of the form
  $f(x)=\F{a}x$ for $a\in\W$.  Then $\Gamma$ satisfies (a) and (b)
above.  Therefore there exists a unique topology on the disjoint union
$\B$ of all $B_x$'s making it into a Banach bundle and such that
  $
  x \mapsto \F {a}x
  $
  is a continuous section for all $a\in\W$.

  \proof
  Property (b) is immediate from the definition of $B_x$.  With
respect to (a) let $f$ be defined by $f(x) = \F{a}x$, where $a\in \W$.
Then
  $$
  \[f(x)\]^2 =
  \[f(x)^* f(x)\] =
  \[\F {a^*}{x\inv} \F ax\].
  $$
  We claim that the map
  $$
  x \mapsto \F{a^*}{x\inv} \F{a}x
  $$
  is continuous as an $\Mult(A)$--valued function on $\Gdual$.  This
will imply that $\[f(x)\]$ depends continuously on $x$.

In order to prove our claim observe that $a^*\rc a$ because $\W$ is
self-adjoint and relatively continuous.  So, for every $x_0$ in
$\Gdual$ we have, with $z=x_0 x\inv$, that
  $$
  \lim_{x\to x_0} \[\F {a^*}{x\inv} \F ax -
    \F {a^*}{x_0\inv} \F a{x_0}\]=
  \lim_{z\to e} \[\F {a^*}{x_0\inv z} \F ax -
    \F {a^*}{x_0\inv} \F a{zx}\] = 0.
  $$
  The last sentence in the statement follows readily from the already
mentioned Theorem 13.18 of \cite{\FD}.
  \proofend

  \state Proposition
  With the multiplication operation $B_x\x B_y \to B_{xy}$ induced by
the corresponding operation on $\Mult(A)$, $\B$ is a Banach algebraic
bundle (as defined in \scite{\FD}{VIII.2.2}).

  \proof
  The only non-trivial axiom left to be verified is the continuity of
the multiplication with respect to the bundle topology.  In order to
verify it we use \scite{\FD}{VIII.2.4}.  That is, it suffices to show
that, given sections
  $$
  \beta(x) = \F ax
  \and
  \gamma(x) = \F bx
  \for x\in\Gdual,
  $$
  with
  $a, b \in \W$, one has that the map
  $$
  (x,y)\in \Gdual\x\Gdual \mapsto \beta(x)\gamma(y) \in \B
  $$
  is continuous.  In order to prove continuity at a given
$(x_0,y_0)\in \Gdual\x\Gdual$ we will use \scite{\FD}{II.13.12} and
hence all we must do is show that there exists a continuous section
$\delta$ such that $\delta(x_0y_0)=\beta(x_0)\gamma(y_0)$ and
  $$
  \lim_{(x,y)\to(x_0,y_0)}
  \[ \beta(x)\gamma(y) - \delta(xy) \] =0.
  $$
  We take $\delta$ to be the section
  $$
  \delta(x) =
  \F{a\F b{y_0}}{x} \for x\in\Gdual.
  $$
  Observing that $a\F b{y_0}$ belongs to $\W$ by spectral invariance
we have that $\delta$ is continuous.
  We have, using \lcite{\FormulasForFourier}.(iv), that
  $$
  \[ \beta(x)\gamma(y) - \delta(xy) \] =
  \[ \F ax \F by - \F{a\F b{y_0}}{xy} \] \$=
  \[\F ax \F b{yy_0\inv y_0} - \F a{xyy_0\inv}\F b{y_0} \],
  $$
  which tends to zero as $z:= yy_0\inv \to e$ (and hence also as
$(x,y)\to(x_0,y_0)$) because $a\rc b$.  This concludes the proof.
  \proofend

It is clear that the adjoint operation is continuous for the bundle
topology by \scite{\FD}{VIII.3.2.(vi$'$)}, and hence $\B$ is in fact a
Banach *-algebraic bundle.  Given that the \cstar-identity ($\[b^*b\]
= \[b\]^2$) obviously holds on $\B$ we conclude that:

  \state Proposition
  $\B$ is a Fell bundle (also referred to as a \cstar-algebraic bundle
in \scite{\FD}{VIII.16.2}).

  This concludes one of the major steps in the development of our
spectral theory which is the description of the spectral subspaces and
of the global topology represented by the Fell bundle $\B$.  The next
step is to relate the cross-sectional \cstar-algebra of $\B$ back to
$A$.

  \state Proposition
  The inclusion map
  $
  \kappa : \B\mapsto \B(\H)
  $
  is a *-representation of $\B$ (in the sense defined in
\scite{\FD}{VIII.9.1}).

  \proof
  The continuity in \scite{\FD}{VIII.8.2.(iv)} is the only point which
we still need to verify.  In order to do so let $\{u_i\}_i$ be a net
in $\B$ converging to some $u_0\in\B$.  Also let $\xi\in\H$ be a
generic vector in $\H$. In order to prove the continuity of our
representation we must therefore prove that
  $u_i\xi\to u_0\xi$ in the norm topology of $\H$.

Let $x_i$ be such that $u_i\in B_{x_i}$.  Given $\varepsilon>0$,
choose a continuous section of\/ $\B$ of the form
  $f(x) = \F ax$, with $a\in\W$, such that
$\[f(x_0)-u_0\]<\varepsilon$.  We have
  $$
  \[u_i\xi - u_0\xi\] \leq
  \[u_i\xi - f(x_i)\xi\] + \[f(x_i)\xi - f(x_0)\xi\] + \[f(x_0)\xi -
u_0\xi\] \$\leq
  \[u_i-f(x_i)\]\[\xi\] + \[f(x_i)\xi - f(x_0)\xi\] +
\varepsilon\[\xi\].
  $$
  Since $f$ is continuous, and so is the norm on $\B$, we have that
  $\[u_i-f(x_i)\]\to \[u_0-f(x_0)\] <\varepsilon$.  The proof will
then be concluded once we show that $\[f(x_i)\xi - f(x_0)\xi\]\to 0$.
  Write $\xi=b\xi'$ for $b\in A$ and $\xi'\in\H$ by the Cohen-Hewitt
factorization Theorem \scite{\HR}{32.22} and observe that
  $$
  \[f(x_i)\xi - f(x_0)\xi\] =
  \[\F a{x_i}b\xi' - \F a{x_0}b\xi'\] \leq
  \[\F a{x_i}b - \F a{x_0}b\]\[\xi'\],
  $$
  which converges to zero by \lcite{\FourierIsStrictlyContinuous}.
  \proofend

  We now wish to show that a certain set of sections is dense in
$L_1(\B)$.  In order to do so we need a result from Classical Harmonic
Analysis for which we have found no specific reference and hence we
prove it below:

  \state Lemma
  \label \HarmonicLemma
  Let $\J(G) = \{g\in C_c(G) : \^g\in L_1(\Gdual)\}.$ Then $\^{\J(G)}$
is dense in $C_0(\Gdual)$.

  \proof
  Observe that $\J(G)$ is a *-subalgebra of $L_1(G)$ (under
convolution).
  Therefore $\^{\J(G)}$ is a *-subalgebra of $C_0(\Gdual)$ (under
pointwise multiplication).
  We claim that $\^{\J(G)}$ satisfies the hypothesis of the
Stone--Weierstrass theorem.
  In fact, suppose that $x_1,x_2\in\Gdual$ are such that $x_1\neq
x_2$.
  Let $t_0\in G$ be such that
  $\du {x_1}{t_0} \neq \du {x_2}{t_0}$
  and take closed convex cones $C_1$ and $C_2$ in $\C$ such that
  $\idu{x_i}{t_0}$ is in the interior of $C_i$ for $i=1,2$, and also
such that $C_1\cap C_2 = \{0\}$.
  Let $\Omega$ be a relatively compact open neighborhood of $t_0$ in
$G$ such that
  $$
  t\in \Omega \quad\Rightarrow\quad \idu {x_i}t\in C_i,
  $$
  for $i=1,2$.
  Using \scite{\HR}{39.16} choose $\psi\in L_1(\Gdual)$ such that
$\^\psi(t_0)=1$ and $\^\psi(G\setminus\Omega) =\{0\}$.  Replacing
$\psi$ by $\psi^*\conv\psi$ we may suppose that $\^\psi(t)\geq 0$ for
all $t$ in $G$.  Let $g = \^\psi$.  Then it is clear that $g\in\J(G)$
and moreover we have
  $$
  \^g(x_i) = \int_G \idu {x_i}t g(t) \d t =
  \int_{\Omega} \idu{x_i}t g(t) \d t \quad\in\quad C_i.
  $$
  Since the integrand is continuous and lies in the interior of $C_i$
for $t=t_0$, it also follows that $\^g(x_i) \neq 0$.  Therefore
$\^g(x_1) \neq \^g(x_2)$, concluding the proof.
  \proofend

  \state Proposition
  \label \DenseSections
  Let $\J(G)$ be as in \lcite{\HarmonicLemma} and let ${\cal S}$ be
the linear span of the set of sections $\eta$ of $\B$ of the form
  $$
  \eta(x) = \^g(x) \F {a}x
  \for x\in \Gdual,
  $$
  where $a\in\W$ and $g\in \J(G)$.  Then $\cal S$ is a dense subset of
$L_1(\B)$ and hence also of $C^*(\B)$.

  \proof
  By \lcite{\TheBanachBundle} and \scite{\FD}{II.13.14} we have that
each $\eta$ of the above form is a continuous section of $\B$.  Since
  $\[\F {a} x\] \leq \[a\]_1$ and $\^g\in L_1(\Gdual)$ it follows that
$\eta\in L_1(\B)$.  Thus ${\cal S} \subseteq L_1(\B)$.

Given $h\in \J(G)$ and $\eta$ as above we have
  $$
  \^h(x)\eta(x) = \^h(x)\^g(x) \F ax =
  \^{(h\conv g)}(x) \F ax
  \for x\in\Gdual.
  $$
  Since $\J(G)$ is an algebra under convolution we conclude that
$\^h\eta$ belongs to $\cal S$.  This says that $\cal S$ is invariant
under pointwise multiplication by $\^h$ for every $h\in \J(G)$.  In
other words
  $\^{\J(G)}{\cal S}\subseteq{\cal S}$.

  Since $\^{\J(G)}$ is dense in $C_0(\Gdual)$ by
\lcite{\HarmonicLemma}, we get $C_0(\Gdual) {\cal
S}\subseteq\overline{\cal S}$, where the closure is taken in
$L_1(\B)$.  One may now easily prove that $\overline{\cal S} \cap
C_c(\B)$ satisfies the hypotheses of \scite{\FD}{II.15.10}, thus
concluding the proof.
  \proofend

The integrated form of $\kappa$ \scite{\FD}{VIII.11.6} is the
representation of the Banach *-algebra $L_1(\B)$ on $\H$, which we
also denote by $\kappa$, given by
  $$
  \kappa(f) v = \int_\Gdual f(x) v \d x
  \for f\in L_1(\B),\quad v\in\H.
  $$
  Since $C^*(\B)$ is defined \scite{\FD}{VIII.17.2} to be the
enveloping \cstar-algebra of $L_1(\B)$ we have that $\kappa$ extends
to $C^*(\B)$.

  \state Proposition
  \label \RangeAndCovariance
  $\kappa(C^*(\B))$ coincides with $\Alg\W$ (Definition
\lcite{\AlgSet}).  Moreover, viewing $\kappa$ as a *-homomorphism from
$C^*(\B)$ to $A$, we have that $\kappa$ is covariant with respect to
the dual action of $G$ on $C^*(\B)$, henceforth denoted by $\beta$,
and the given action $\a$ on $A$.

  \proof
  Let $a\in\W$ and let $g\in \J(G)$.  We then have that the element
  $$
  a':=\int_G g(t) \a_t(a) \d t
  $$
  clearly belongs to $\Alg\W$ and is $\a$-integrable by
\lcite{\IntIsIntegrable}.  Observe that the section
  $$
  \eta : x\in \Gdual \longmapsto \F {a'}x
\={\FourierTransfOfConvolution} \^g(x) \F ax \in B_x
  $$
  is in the space $\cal S$ introduced in \lcite{\DenseSections}.  In
fact $\cal S$ is spanned by the section of this form.  Since $\eta\in
L_1(\B)$ we have by \lcite{\InvFourierTransform} that
  $$
  a' =
  \int_\Gdual \F {a'}x \d x =
  \int_\Gdual \eta(x) \d x =
  \kappa(\eta).
  $$
  This shows that $a'\in\kappa(C^*(\B))$ as well as that $\kappa({\cal
S})\subseteq \Alg\W$.  Using \lcite{\DenseSections} we then conclude
that $\kappa(C^*(\B))\subseteq \Alg\W$.

Applying \scite{\HR}{39.16} (to the dual group $\Gdual$) we may find
$g\in \J(G)$ such that $\[a'-a\]$ is arbitrarily small.  This shows
that $a$ lies in the range of $\kappa$ and hence that $\W\subseteq
\kappa (C^*(\B))$.  Summarizing our findings so far we have
  $$
  \W\subseteq \kappa (C^*(\B))\subseteq \Alg\W.
  \leqno{(\seqnumbering)}
  \label \TwoInclusions
  $$

Let $a\in\W$, $g\in \J(G)$, and consider $a'$ and $\eta$ defined in
terms of $a$ and $g$ as above.  We claim that $\a_t(\kappa(\eta)) =
\kappa(\beta_t(\eta))$ for all $t\in G$, where $\beta$ refers to the
dual action.  In order to prove it observe that
  $$
  \kappa(\beta_t(\eta)) =
  \int_\Gdual \beta_t(\eta)\calcat x \d x =
  \int_\Gdual \idu xt \eta(x) \d x \$=
  \int_\Gdual \idu xt \F {a'}x \d x \={\InvFourierTransform}
  \a_t(a') = \a_t(\kappa(\eta)).
  $$
  This proves our claim and, since $\cal S$ is dense in $C^*(\B)$, it
also proves the last assertion in the statement.  In addition we
conclude that the range of $\kappa$ is invariant under $\a$ and hence,
in view of \lcite{\TwoInclusions}, we have that $\kappa (C^*(\B)) =
\Alg\W$.
  \proofend

  Consider the representation $\kappa\*\lambda$ of $\B$ on
$L_2(\Gdual,\H)$.  This is given by
  $$
  (\kappa\*\lambda)(b_x) = \kappa(b_x)\*\lambda_x
  \for x\in \Gdual,\quad b_x\in B_x,
  $$
  where $\lambda$ refers to the left regular representation of
$\Gdual$.  As before we also denote by $\kappa\*\lambda$ the
corresponding integrated representation of $C^*(\B)$.

Recall from section \lcite{\BasicFactsSection} that $\pi$ is the
representation of $A$ on $L_2(G,\H)$ given by
  $
  \pi(a)\xi\calcat t = \a_t\inv (a)\xi(t),
  $
  for $a\in A$, $\xi\in L_2(G,\H)$, and $t\in G$.

  \state Proposition
  \label \PropositionOfDiagram
  Let $a':=\int_G g(t) \a_t(a) \d t$, where $a\in \W$ and $g\in
\J(G)$, and let $\eta$ be the element of $C^*(\B)$ represented by the
section
  $\eta(x)
  = \^g(x) \F ax
  = \F {a'}x$ for $x\in\Gdual$.  Then the following diagram commutes
  $$
  \matrix{
  &\pi(a')\cr
  L_2(G,\H) & \bigarrow & L_2(G,\H) \cr\cr
  \Ft \Big\downarrow && \Big\downarrow \Ft\cr\cr
  L_2(\Gdual,\H) & \bigarrow & L_2(\Gdual,\H) \cr
  &{(\kappa\*\lambda)\eta}\cr
  }
  $$
  where $\Ft$ stands for the Fourier Transform.

  \proof
  Let $\xi\in C_c(G,\H)$ and $x\in\Gdual$.  Then
  $$
  \Ft \pi(a') \xi \calcat x=
  \int_G \idu xt \pi(a') \xi \calcat t \d t =
  \int_G \idu xt \a_t\inv(a') \xi(t) \d t =
  $$
  $$
  \={\InvFourierTransform}
  \int_G \idu xt \(\int_\Gdual \idu y{t\inv} \F {a'}y \d y\) \xi(t) \d
t =
  \int_\Gdual \F {a'}y \int_G \idu {y\inv x}t \xi(t) \d t\d y \$=
  \int_\Gdual \F {a'}y \^\xi(y\inv x) \d y =
  \int_\Gdual \(\F {a'}y \* \lambda_y\) \^\xi\,\calcat x \d y =
  ((\kappa\*\lambda)\eta) \^\xi\,\calcat x.
  $$
  Since $x$ is arbitrary we have that
  $ \Ft \pi({a'}) \xi = ((\kappa\*\lambda)\eta)\Ft \xi$.  Since
$C_c(G,\H)$ is dense in $L_2(G,\H)$ the proof is concluded.
  \proofend

  \state Proposition
  The representation $\kappa$ of $C^*(\B)$ defined above is faithful
and hence $C^*(\B)$ is covariantly isomorphic to $\Alg\W$.

  \proof
  Let $a':=\int_G g(t) \a_t(a) \d t$, where $a\in \W$ and $g\in
\J(G)$.  Set
  $\eta(x)
  = \^g(x) \F ax
  = \F {a'}x$ for $x\in\Gdual$,
  so that $\eta\in \cal S$ and $\kappa(\eta) = a'$, as before.
Consider the composition of maps
  $$
  C^*(\B) \buildrel \kappa \over \longrightarrow
  A\buildrel \pi \over \longrightarrow
  \B(L_2(G,\H)) \buildrel {{\scriptstyle Ad}_\Ft}\over \longrightarrow
  \B(L_2(\Gdual,\H)),
  $$
  where $Ad_\Ft$ is the conjugation by the Fourier transform.
  We then have that
  $$
  Ad_\Ft(\pi(\kappa(\eta))) =
  \Ft \pi(a')\Ft^* =
  (\kappa\*\lambda)(\eta),
  $$
  where the last step follows from \lcite{\PropositionOfDiagram}.
  Since the set of $\eta$'s considered span $\cal S$, and since $\cal
S$ is dense in $C^*(\B)$, we conclude that the above composition of
maps coincides with $\kappa\*\lambda$, which is a faithful
representation of $C^*(\B)$ by
  \scite{\Rep}{3.6}.
  This shows that $\kappa$ is one to one and hence the proof is
complete.
  \proofend

  The following Theorem subsumes our findings in this section and is
one of our main results.  We observe that whereas we have worked above
under the assumption that $\W$ is a spectrally invariant linear space,
these hypotheses are not needed below.

  \state Theorem
  \label \SpectralTheorem
  Let $\a$ be a strongly continuous action of a locally compact
abelian group $G$ on a \cstar-algebra $A$ and let $\W$ be a subset of
$A$ such that
  \izitem $\W^*=\W$,
  \zitem $\W$ consists of $\a$-integrable elements (Definition
\lcite{\AltDefinitionOfIntegr}) and
  \zitem $\W$ is relatively continuous (Definition
\lcite{\RelativelyContinuousSpectrallyInvariant}).
  \medskip\noindent
  Then there exists a Fell bundle $\B$ over $\Gdual$ such that
$C^*(\B)$ is isomorphic to the smallest $\a$-invariant
sub-\cstar-algebra of $A$ containing $\W$ under an isomorphism which
is covariant with respect to $\a$ and the dual action of $G$ on
$C^*(\B)$.

  \proof
  Let $\widetilde\W$ be the spectrally invariant hull of $\W$.
Observing the description of $\widetilde\W$ given in the proof of
\lcite{\ConsertaMasNaoEstraga} it is clear that $\widetilde\W$ is
self-adjoint.  Also by \lcite{\ConsertaMasNaoEstraga} we have that
$\widetilde\W$ is relatively continuous, and hence $\widetilde\W$
satisfies \lcite{\WorkingHypothesis}.(i--iv).  It follows that
$\span(\widetilde\W)$ satisfies all of the conditions in
\lcite{\WorkingHypothesis}.

Note that, by \lcite{\SIHHasSameAlgebra} we have that
  $\Alg{\W} = \Alg{\widetilde\W} = \Alg{\span(\widetilde\W)}$.  The
conclusion then follows from the previous results applied to
$\span(\widetilde\W)$.
  \proofend

  Combining the previous result with \lcite{\RelContInBundles} we
obtain the following:

  \state Corollary
  \label \CharacterizationOfDualaction
  Let $(A,G,\a)$ be a \cstar-dynamical system where $G$ is abelian.
  A necessary and sufficient condition for $\a$ to be equivalent to a
dual action is that $A$ contains a dense, self-adjoint, relatively
continuous set of $\a$-integrable elements.

  We have not been able to determine the extent to which condition
(iii) in \lcite{\SpectralTheorem} is really necessary and hence we
leave open the following:

  \sysstate{Question}{\rm}
  {\label \MainQuestion Suppose that the set of $\a$-integrable
elements is dense in $A$, or even that the set of linear combinations
of positive $\a$-integrable elements is
dense\footnote{\footcntr}{\eightpoint This is Rieffel's tentative
definition of proper actions \scite{\RieffelIntegrable}{4.5}.}.  Does
if follow that there exists a dense subset $\W$ of $\a$-integrable
elements which is relatively continuous?}

  \section {Classical dynamical systems}
  \label \ClassicalDynSystemsSection
  In this section we shall apply the results obtained so far to the
case of a classical dynamical system, that is, a group action on a
locally compact topological space.  Given that our methods were
developed for abelian groups we shall let, throughout this section,
$G$ be a locally compact abelian group and
  $$
  \a : (t,p)\in G\x \T\quad \longmapsto\quad tp\in \T
  $$
  be a continuous action of $G$ on a locally compact space $\T$.  One
therefore gets a strongly continuous action of $G$ on the
\cstar-algebra
  $A=C_0(\T)$, which we will also denote by $\a$, given by
  $$
  \a_t(f)\calcat p = f(t\inv p)
  \for t\in G,\quad f\in C_0(\T),\quad p\in \T.
  $$

  Let $f\in C_0(\T)$ be an $\a$-integrable element.  Then, seeing
$\Mult(C_0(\T))$ as the set $C_b(\T)$ of bounded continuous functions
on $\T$, we have that $\F fe$ is represented by the function
  $$
  \F fe(p) = \int_G f(t\inv p) \d t
  \for p\in\T.
  $$
  By a simple change of variables we see that $\F fe$ is constant on
the orbits of $\T$ under $G$ and hence it defines, by passage to the
quotient, a continuous function on the orbit space $\T/G$, even if
this not a locally compact or Hausdorff space!  In any case $\T/G$
carries the quotient topology and the
  the quotient map
  $$
  Q: \T \to \T/G
  $$
  is continuous.

  \definition
  We shall denote by $C_0(\T/G)$ the space of continuous complex
valued functions $f$ on $\T/G$ such that for every $\varepsilon>0$,
there exists a compact subset $K\subseteq\T$ such that
$|f(q)|<\varepsilon$, for all
  $q\in( \T/G) \setminus Q(K)$.

  That is, $C_0(\T/G)$ is defined in the usual way, except that when
it comes to considering compact subsets of $\T/G$ we take only the
images of compact subsets of $\T$ under the quotient map.
Equivalently, $C_0(\T/G)$ consists of the continuous functions $f$ on
$\T$ which are constant along every orbit and such that for every
$\varepsilon>0$, there exists a compact subset
  $K\subseteq\T$ such that
  $|f(p)|<\varepsilon$, for every $p\in\T$ outside the orbit $\O(K)$
of $K$.

  \state Lemma
  \label \RelContForFunctions
  Let $f$ be an $\a$-integrable element of $C_0(\T)$.  Suppose that
$f$ is positive and that $\F fe$ belongs to $C_0(\T/G)$. Then $f$ is
\stress{absolutely continuous} in the sense that $f\rc g$ for every
$\a$-integrable $g$.

  \proof Fix an $\a$-integrable element $g\in C_0(\T)$.
  Let $\varepsilon>0$ be given and take a compact subset
$K\subseteq\T$ such that $|\F fe(p)|<\varepsilon/(2\[g\]_1)$, for
every $p\in\T\setminus\O(K)$.
  Observe that for all $p\in \T$ and $w\in\Gdual$ we have
  $$
  |\F fw(p)| \leq
  \int_G \left|\du wt f(t\inv p)\right| \d t =
  \int_G f(t\inv p) \d t =
  \F fe(p).
  $$
  Thus $|\F fw(p)| < \varepsilon/(2\[g\]_1)$, for all
$p\in\T\setminus\O(K)$.
  For any such $p$, and any $x,y,z\in\Gdual$, we therefore have that
  $$
  \left|\big(\F{f}{xz}\F{g}{y} - \F{f}{x}\F{g}{zy}\big)\calcat
p\right| \leq
  |\F{f}{xz}(p)|\[g\]_1 + |\F{f}{x}(p)| \[g\]_1 \leq
  \varepsilon.
  $$

  Recall that for every $w\in\Gdual$, one has that $\F
fw\in\Mult_w(C_0(\T))$ and hence for any $p\in\T$ and $s\in G$,
  $
  \F fw(s\inv p) =
  \idu ws \F fw(p).
  $
  It follows that
  $$
  \big(\F{f}{xz}\F{g}{y} - \F{f}{x}\F{g}{zy}\big)\calcat {s\inv p} =
  \idu{xzy}s \big(\F{f}{xz}\F{g}{y} - \F{f}{x}\F{g}{zy}\big)\calcat p,
  $$
  for every $x,y,z\in\Gdual$.
  For brevity we will use the abbreviation
  $$
  \Dxyz = \F{f}{xz}\F{g}{y} - \F{f}{x}\F{g}{zy}.
  $$
  We have therefore proven that, for every $x,y,z\in\Gdual$
  \izitem if $p\notin\O(K)$ then $|\Dxyz(p)|\leq \varepsilon$, and
  \zitem if $p\in\T$ and $s\in G$ then $\Dxyz(s\inv p) = \idu{xzy}s
\Dxyz(p)$.
  \medskip\noindent
  Let $h\in C_0(\T)$ be such that $h(p)=1$ for every $p\in K$, so
that, using (ii),
  $$
  \sup_{p\in\O(K)} |\Dxyz(p)| =
  \sup_{p\in K} |\Dxyz(p)| \leq
  \[\Dxyz h\] =
  \[\F{f}{xz}\F{g}{y}h - \F{f}{x}\F{g}{zy}h\] \$\leq
  \[\F{f}{xz}\F{g}{y}h - \F{f}{x}\F{g}{y}h\] +
  \[\F{f}{x}\F{g}{y}h - \F{f}{x}\F{g}{zy}h\] \$\leq
  \[g\]_1\[\F{f}{xz}h - \F{f}{x}h\] +
  \[f\]_1\[\F{g}{y}h - \F{g}{zy}h\].
  $$

  By \lcite{\FourierIsStrictlyContinuous} there exists a neighborhood
$\Omega$ of $e$ in $\Gdual$ such that the above is less than
$\varepsilon$ for all $z\in\Omega$ and all $x,y\in\Gdual$.  We
conclude that, for $z\in\Omega$,
  $$
  \sup_{x,y\in\Gdual}
  \[\F{f}{xz}\F{g}{y} - \F{f}{x}\F{g}{zy}\] =
  \sup_{\buildrel{\scriptstyle p\in\T}\over
    {\scriptstyle x,y\in\Gdual}} |\Dxyz(p)| \leq \varepsilon,
  $$
  and hence that $f\rc g$.
  \proofend

Recall from \scite{\RieffelIntegrable}{4.7} that $\a$ is a proper
action in the usual sense\footnote{\footcntr}
  {\eightpoint A group action is said to be proper when the map
$(t,x)\in G\x \T \mapsto (tx,x)\in \T\x \T$ is proper, i.e the inverse
image of compact sets is compact.}
  if and only if the linear span of the set of positive
$\a$-integrable elements is dense in $C_0(\T)$ (see remark
\lcite{\RemarkOnDifferentTerminology} regarding a difference between
our terminology and that adopted in \cite{\RieffelIntegrable}).  In
fact, for a proper action any $f$ in $C_c(\T)$ is $\a$-integrable (see
\scite{\RieffelProper}{Section 1}).

  \state Proposition
  \label \AbsCont
  If $\a$ is proper then any $f$ in $C_c(\T)$ is absolutely continuous
(as defined in \lcite{\RelContForFunctions}).

  \proof
  Let $f\in C_c(\T)$.  It is then clear that
  $\F fe$ is compactly supported in $\T/G$ and, in particular, that
  $\F fe\in C_0(\T/G)$.  Therefore if $f$ is also positive the
conclusion follows from \lcite{\RelContForFunctions}.  In the general
case observe that $f$ is a linear combination of four positive
elements in $C_c(\T)$ each of which is absolutely continuous.
  \proofend

Let $\P = \A_+ \cap C_c(\T)$.  Then by the above result $\P$ is a
relatively continuous set.  It is also clear that $\P$ is a hereditary
cone and satisfies \lcite{\PropertiesOfMaximalCones}.(ii).  We may
then use \lcite{\AbstractMR} to derive the following well known result
(see \cite{\Green}, \cite{\RieffelTransf}, \cite{\RieffelIntegrable}):

  \state Corollary
  Let $\a$ be a proper action of the abelian group $G$ on the locally
compact topological space $\T$.  Then the closure of $\zeta(C_c(\T))$
is the imprimitivity bimodule for a {\MR} equivalence between
$C_0(\T/G)$ and an ideal in $\A\crossproduct_\a G$.

  \proof
  Let $\P = \A_+ \cap C_c(\T)$ as above.  Then it is clear that the
sets $\N$ and $\M$ constructed in terms of $\P$ as in
\lcite{\AbstractMR} both coincide wih $C_c(\T)$.
  By \lcite{\AbstractMR} all we need to verify is that the closure of
$\F {C_c(\T)}e$ coincides with $C_0(\T/G)$.  As observed in the proof
of \lcite{\AbsCont} we have that $\F {C_c(\T)}e\subseteq C_c(\T/G)$.
  The conclusion then follows easily from the Stone--Weierstrass
Theorem.
  \proofend

Let us now consider applying Theorem \lcite{\SpectralTheorem} to the
present situation.  For this, observe that taking $\W=C_c(\T)$, all of
the hypotheses of \lcite{\SpectralTheorem} hold.  Moreover, since $\W$
is dense in $C_0(\T)$, we have that $\Alg\W=C_0(\T)$.  Theorem
\lcite{\SpectralTheorem} then immediately implies:

  \state Corollary
  \label \SpectralTheoremCommutativeCase
  Let $\a$ be a proper action of a locally compact abelian group $G$
on a locally compact topological space $\T$.  Then $\a$ is equivalent
to a dual action.  That is, there exists a Fell bundle $\B$ over
$\Gdual$ such that $C_0(\T)$ is covariantly isomorphic to
  $C^*(\B)$ with respect to the dual action of $G$ on $C^*(\B)$.

  \section {An example}
  Throughout this article we have been working under the tacit
assumption that dual actions are the ``good guys'' among
\cstar-dynamical systems, and we hope to have given enough evidence to
convince the reader of this.  Nevertheless we would now like to
explore a specific example of dual action which will bring to light
some of the complexities that may be found as well as some of the
inherent difficulties in answering the questions that we have posed.

Our example will in fact be that of a classical dual action as
described in \scite{\Ped}{7.8.3}.  Consider the action of the unit
circle $S^1$ on $C(S^1)$ by translation.  The crossed product algebra
is well known to be isomorphic to the algebra $\K$ of compact
operators on $\ell_2(\Z)$.  It is also well known that the dual action
of $\Z$ on $\K$ is given by conjugation by the powers of the forward
bilateral shift (this is in fact a consequence of Takai duality since
the first action described above is already a dual action with respect
to the trivial action of $\Z$ on a point).

The action we wish to discuss is thus the action of $\Z$ on $\K$ given
by
  $$
  \a_n(T) = U^nTU^{-n}
  \for n\in\Z,\quad T\in\K,
  $$
  where $U$ is the forward bilateral shift on $\ell_2(\Z)$.

  In order to avoid an excessive use of the Fourier transform we will
identify the Hilbert spaces $\ell_2(\Z)$ and $L_2(S^1)$, henceforth
denoted simply by $\H$, by considering the usual orthonormal basis
$\{\e n\}_{n\in\Z}$ of the latter given by $\e n(z) = z^n$ for all $z$
in $S^1$.  Under this identification $U$ becomes the operator of
pointwise multiplication by $z$:
  $$
  U\xi\calcat z = z\xi(z)
  \for\xi\in L_2(S^1),\quad z\in S^1.
  $$

  \state Proposition
  \label \RankOneProjection
  Let $P$ be a rank-one projection on $\H$.  Then $P$ is
$\a$-integrable if and only $P$ has the form
  $$
  P(\xi) = \<\xi,\phi\>\phi
  \for \xi\in\H,
  $$
  where $\[\phi\]_2=1$ and $\phi$ belongs to $L_\infty(S^1)$, viewed
as a subspace of $L_2(S^1)$.  In that case the series
  $
  \sum_{n\in\Z} \a_n(P)
  $
  converges strictly-unconditionally to the operator of pointwise
multiplication by $|\phi|^2$.

  \proof
  Suppose that $P$ has the form mentioned in the statement for some
$\phi\in L_\infty(S^1)$ and denote by $M_\phi$ the operator on $\H$
given by pointwise multiplication by $\phi$.
  Also, given any finite subset $J\subseteq\Z$, let $Q_J$ denote the
orthogonal projection onto $\span\{\e n \colon n\in J\}$.  We then
have for all $\xi\in\H$ that
  $$
  M_\phi Q_J M_\phi^* (\xi) =
  M_\phi \( \sum_{n\in J} \<M_\phi^* \xi,\e n\> \e n \) =
  \sum_{n\in J} \<\xi,\phi \e n\> \phi \e n \$=
  \sum_{n\in J} \<\xi,U^n(\phi)\> U^n(\phi) =
  \sum_{n\in J} U^n P U^{-n} (\xi) =
  \sum_{n\in J} \a_n(P) (\xi).
  $$
  It follows that
  $
  M_\phi Q_J M_\phi^* = \sum_{n\in J} \a_n(P).
  $

  If $T$ is any compact operator on $\H$ we have that
  $\lim_J Q_JT = \lim_J TQ_J = T$, where the limit is with respect to
the directed set formed by all compact, i.e. finite, subsets
$J\subseteq \Z$.  Therefore, still assuming that $T$ is compact,
  $$
  \lim_J \sum_{n\in J} \a_n(P) T =
  \lim_J M_\phi Q_J M_\phi^* T =
  M_\phi M_\phi^* T =
  M_{|\phi|^2} T,
  $$
  and similarly if $T$ is on the left hand side.  It follows that the
series
  $
  \sum_{n\in \Z} \a_n(P)
  $
  converges strictly-unconditionally to
  $M_{|\phi|^2}$ proving that $P$ is $\a$-integrable, and also proving
the last sentence in the statement.

  Conversely suppose that the rank-one projection $P$ is
$\a$-integrable and let $\phi$ be a vector in the range of $P$ of norm
one.  Thus
  $
  P(\xi) = \<\xi,\phi\>\phi,
  $
  for $\xi\in\H$.
  Given $i,j$ in $\Z$ we have
  $$
  \<P \e j, \e i\> =
  \<\<\e j,\phi\> \phi,\e i\> =
  \^\phi(i) \overline{\^\phi(j)} =
  \^\phi(i) \^{\overline\phi}(-j),
  $$
  where $\overline\phi$ is the pointwise complex conjugate of $\phi$.
  Let $T = \sum_{k\in \Z} \a_k(P)$, where the series converges
strictly-unconditionally and hence also in the weak-operator topology.
We then have
  $$
  \<T \e j, \e i\> =
  \sum_{k\in\Z} \<\a_k(P) \e j, \e i\> =
  \sum_{k\in\Z} \<U^k P U^{-k} \e j, \e i\> =
  \sum_{k\in\Z} \<P \e {j-k}, \e {i-k}\> \$=
  \sum_{k\in\Z} \^\phi(i-k) \^{\overline\phi}(k-j) =
  \(\^\phi * \^{\overline\phi}\) (i-j) =
  \^{\phi\overline\phi}(i-j).
  $$
  It follows that the matrix of $T$ with respect to the canonical
basis $\{\e n\}_{n\in\Z}$ is a Laurent matrix (in the classical sense
\scite{\Halmos}{241}) with symbol $|\phi|^2$.  Therefore, since $T$ is
bounded, we must have that $\phi\in L_\infty(S^1)$.
  \proofend

  \sysstate{Remark}{\rm}
  {The first important consequence to be drawn from this result is
that Definition 6.4 in \cite{\RieffelIntegrable} leads to a
\stress{{\GFPA}} bigger than desired.  In fact, adopting that
definition, by \lcite{\RankOneProjection} we would have that the
{\GFPA} will be formed by all multiplication operators by measurable
bounded functions on $S^1$ and hence will be isomorphic to
$L_\infty(S^1)$.  However, since we are in a situation where the
classical Takai duality Theorem holds, there seems to be little doubt
that the ``correct'' definition of {\GFPA} should yield $C(S^1)$
instead.  In addition our example illustrates that taking the
integrals of higher powers of integrable elements is not the correct
approach either since the trouble here is caused by idempotent
operators!  We should note that this was also indicated in
\scite{\unconditional}{Section 7}.}

As a simple consequence of \lcite{\RankOneProjection} we can give a
precise characterization of positive integrable elements:

  \state Proposition
  \label \PositiveOperators
  Let $T$ be a positive compact operator on $\H$. Then $T$ is
$\a$-integrable if and only if there exists a sequence
$\{\lambda_n\}_n$ of positive real numbers converging to zero and a
pairwise orthogonal sequence $\{\phi_n\}_n \subseteq L_\infty(S^1)$
such that
 $\sum_{n=1}^\infty \lambda_n |\phi_n|^2$ is pseudo-summable in
$L_\infty(S^1)$ (meaning that the finite sums are uniformly bounded
{\rm \scite{\unconditional}{2.5}}), and
  $$
  T = \sum_n \lambda_n P_n,
  $$
  where $P_n$ is the projection onto the subspace of $\H$ spanned by
$\phi_n$.

  \proof
  Suppose that $T$ is $\a$-integrable.  Using the spectral Theorem for
compact self-adjoint operators write $T = \sum_n \lambda_n P_n$, where
$\lambda_n>0$ for all $n$, $\lim_n\lambda_n=0$, and the $P_n$ are
pairwise orthogonal rank-one projections.  Since the set of
$\a$-integrable elements is a hereditary cone by
\lcite{\BasicPropertiesOfModule}, we have that each $P_n$ is
$\a$-integrable and hence, by the above Lemma, there are bounded
measurable functions $\phi_n$, necessarily pairwise orthogonal, such
that $P_n$ is the projection onto the one dimensional space spanned by
$\phi_n$.  For each integer $N$ we have that
  $
  0 \leq \sum_{n=1}^N \lambda_n P_n \leq T,
  $
  and hence
  $$
  0 \leq \E{\sum_{n=1}^N \lambda_n P_n} =
  \sum_{n=1}^N \lambda_n M_{|\phi_n|^2} \leq \E T,
  $$
  which implies that
  $
  \sum_{n=1}^N \lambda_n |\phi_n|^2
  $
  is uniformly bounded with $N$.

  Conversely, suppose that $T$ is of the above form and let us prove
it to be $\a$-integrable.  For a given finite subset $J\subseteq\Z$ we
have
  $$
  \sum_{k\in J} \a_k(T) =
  \sum_{k\in J} \sum_{n\in\Z} \lambda_n \a_k(P_n) =
  \sum_{n\in\Z} \lambda_n \sum_{k\in J} \a_k(P_n) \$=
  \sum_{n\in\Z} \lambda_n M_{\phi_n} Q_J M_{\phi_n}^* \leq
  \sum_{n\in\Z} \lambda_n M_{|\phi_n|^2},
  $$
  where $Q_J$ is as in the proof of \lcite{\RankOneProjection}.  It
follows that the net $\{\sum_{k\in J} \a_k(T)\}_J$ is bounded and
since it is also increasing it must converge strongly.  The strong
topology can be easily shown to coincide with the strict topology on
bounded subsets of $\K$ and hence we have shown that our net converges
strictly.  This completes the proof.
  \proofend

Let us now consider the problem mentioned near the end of Section
\lcite{\LeftStructure} as to whether $\zeta(a)\zeta(b)^*$ belongs to
the crossed product algebra for a given pair of elements $a,b\in\Na$.
We will restrict our attention to rank-one projections since this will
suffice to illustrate a point to be made below and also because the
behavior of general elements of $\K$ is mirrored in the behavior of
rank-one projections, as seen in \lcite{\PositiveOperators}.  We must
first compute the Fourier transform:

  \state Proposition
  \label \ComputingFT
  Let $\phi\in L_\infty(S^1)$ with $\[\phi\]_2=1$, and let $P$ be the
projection onto the subspace of $\H$ spanned by $\phi$.  Then for each
$z$ in the Pontryagin dual\/ $S^1$ of\/ $\Z$ one has
  $$
  \F Pz = M_\phi W_z M_\phi^*,
  $$
  where $W_z$ is the diagonal operator given on the canonical basis by
  $$
  W_z(\e n) = z^{n}\e n \for n\in\Z.
  $$

  \proof
  In order to prove the statement it suffices to show that the
operators above have identical matrices with respect to the canonical
basis $\{\e n\}_{\in\Z}$.  Recaling that
  $\F Pz = \sum_{n\in\Z} z^{n} \a_n(P)$, where the sum converges
strictly-unconditionally, we have
  $$
  \<\F Pz \e j, \e i\> =
  \sum_{n\in\Z} z^{n} \<U^n P U^{-n} \e j, \e i\> =
  \sum_{n\in\Z} z^{n} \<P \e {j-n}, \e {i-n}\> \$=
  \sum_{n\in\Z} z^{n} \<\<\e{j-n},\phi\>\phi, \e {i-n}\> =
  \sum_{n\in\Z} z^{n} \^\phi(i-n) \overline{\^\phi(j-n)}.
  $$
  On the other hand
  $$
  \<M_\phi W_z M_\phi^* \e j, \e i\> =
  \<W_z M_\phi^* \e j, M_\phi^* \e i\> =
  \sum_{n\in\Z} z^{n} \<M_\phi^* \e j,\e n\>
    \overline{\langle M_\phi^* \e i,\e n\rangle} \$=
  \sum_{n\in\Z} z^{n} \<\e j,M_\phi \e n\>
    \langle M_\phi \e n, \e i\rangle =
  \sum_{n\in\Z} z^{n} \overline{\^\phi(j-n)} \^\phi(i-n) ,
  $$
  concluding the proof.
  \proofend

  \state Proposition
  \label \RCForRankOneOperators
  Let $\phi$ and $\psi$ be in $L_\infty(S^1)$ with
$\[\phi\]_2=\[\psi\]_2=1$ and denote by $P$ and $Q$ the projections
onto the one-dimensional subspaces of $\H$ spanned by $\phi$ and
$\psi$ respectively.  Then the following are equivalent
  \izitem $\zeta(P)\zeta(Q)^*$ belongs to $\K\crossproduct_\a \Z$,
  \zitem $P\rc Q$,
  \zitem $\overline\phi\psi\in C(S^1)$.

  \proof
  Given that $P=P^*P$ and similarly for $Q$, the equivalence between
(i) and (ii) follows immediately from
\lcite{\ReformulationOfIncredibleFormula}.  Therefore we need only
prove that (ii) and (iii) are equivalent.
  During the course of this proof we will denote the multiplication
operator $M_f$ simply by $f$ for each $f$ in $L_\infty$.
 Let $\tau$ be the (discontinuous) action of $S^1$ on $L_\infty(S^1)$
given by
  $$
  \tau_z f\calcat w = f(z w)
  \for z,w\in S^1, \quad f\in L_\infty(S^1),
  $$
  and observe that for all $f$ in $L_\infty(S^1)$ one has
  $
  W_z f W_z^* = \tau(f).
  $
  Given $x,y,z\in S^1$ we have
  $$
  \[\F P{xz} \F Qy - \F Px \F Q{zy}\] =
  \[\phi W_{xz} \phi^* \psi W_y\psi^* -
    \phi W_x \phi^* \psi W_{zy} \psi^*\] \$\leq
  \[\phi W_x\]
  \[W_z \phi^*\psi - \phi^*\psi W_z \]
  \[W_y\psi^*\] =
  \[\phi\]_\infty \[\psi\]_\infty
  \[\tau_z(\overline\phi\psi) - \overline\phi\psi \]_\infty.
  $$
  Thus, assuming that $\overline\phi\psi\in C(S^1)$, we have that
$P\rc Q$.
  Conversely, suppose that $P\rc Q$.  Then for all $z$ in $S^1$ we
have
  $$
  \[\tau_z(\overline\phi\psi) - \overline\phi\psi\]^3 =
  \[
  \big(W_z\phi^*\psi W_z^* - \phi^*\psi\big)
  \big(W_z\psi^*\phi W_z^* - \psi^*\phi\big)
  \big(W_z\phi^*\psi W_z^* - \phi^*\psi\big)
  \]
  \$\leq
  \big\|
  W_z\phi^*\psi \psi^*\phi \phi^*\psi W_z^*
    - W_z\phi^*\psi W_z^* \psi^*\phi W_z\phi^*\psi W_z^*
  +$$$$- W_z\phi^*\psi \psi^*\phi W_z^* \phi^*\psi
    + W_z\phi^*\psi W_z^* \psi^*\phi \phi^*\psi
  +$$$$- \phi^*\psi W_z\psi^*\phi \phi^*\psi W_z^*
    + \phi^*\psi \psi^*\phi W_z\phi^*\psi W_z^*
  \$+ \phi^*\psi W_z\psi^*\phi W_z^* \phi^*\psi
    - \phi^*\psi \psi^*\phi \phi^*\psi
  \big\|
  \$\leq
  \[\phi\]_\infty\[\psi\]_\infty\Big(
  \[\psi\psi^*\phi\phi^*
    -\psi W_z^*\psi^*\phi W_z\phi^*\] +
  \[\psi\psi^*\phi W_z^*\phi^*
    -\psi W_z^*\psi^*\phi\phi^*\]\$+
  \[\psi W_z\psi^*\phi\phi^*
  -\psi\psi^*\phi W_z\phi^*\] +
  \[\psi W_z\psi^*\phi W_z^*\phi^*
    -\psi\psi^*\phi\phi^*\]
  \Big)
  \$=
  \[\phi\]_\infty\[\psi\]_\infty\Big(
  \[\E Q \E P
    -\F Q{z\inv}\F Pz\] +
  \[\E Q \F P{z\inv}
    -\F Q{z\inv}\E P \]\$+
  \[\F Qz\E P
    -\E Q \F Pz\] +
  \[\F Qz\F P{z\inv}
    -\E Q \E P \]
  \Big),
  $$
  which tends to zero as $z\to 1$ by (ii).  This shows that
$\overline\phi\psi$ is continuous.
  \proofend

  Given that absolutely continuous elements played an important role
in the case of classical dynamical systems, studied above, it is
interesting to characterize them in the present situation.  However we
have:

  \state Proposition
  \label \NoPositiveAbsolutelyContinuous
  Let $T$ be a positive $\a$-integrable element of\/ $\K$.  Suppose
that $T$ is absolutely continuous in the sense that $T\rc S$ for all
$\a$-integrable $S\in\K$.  Then $T=0$.

  \proof
  Suppose by contradiction that $T\neq0$ and pick a nonzero vector
$\phi$ in $\H$ such that $T(\phi) = \[T\]\phi$.  Also let $P$ be the
orthogonal projection onto the one-dimensional space spanned by
$\phi$, so that $\[T\]P \leq T$.  So $P$ is $\a$-integrable by
\lcite{\BasicPropertiesOfModule}.(i) and hence $\phi\in L_\infty(S^1)$
by \lcite{\RankOneProjection}.

Given $\psi\in L_\infty(S^1)$ let $Q$ be the orthogonal projection
onto the one-dimensional subspace of $\H$ spanned by $\psi$.  Then
  $$\[T\]P \leq T\rc Q$$
  which implies by \lcite{\Herege} that
  $P \rc Q$.  Therefore, using \lcite{\RCForRankOneOperators} we have
that $\overline\phi\psi\in C(S^1)$.  But this can only happen for all
$\psi$ in $L_\infty(S^1)$ if $\phi=0$.
  \proofend

  As already mentioned $\K$ is the crossed product algebra
$C(S^1)\crossproduct_\tau S^1$, where $\tau$ is given by translation.
Therefore $\K$ is also the cross-sectional \cstar-algebra for the
semi-direct product bundle \scite{\FD}{VIII.4} $\B$ whose total space
is
  $C(S^1)\x S^1$, carrying the bundle operations
  $$
  (f,z)\cdot(g,w) = (f\tau_z(g),zw)
  \and
  (f,z)^* = (\tau_z\inv(\overline f),z\inv),
  $$
  for $f,g\in C(S^1)$ and $z,w\in S^1$.  Consider the representation
of $\B$ on $\H=L_2(S^1)$ given by
  $$
  \rho : (f,z)\in\B \mapsto M_f W_z\in \B(\H),
  $$
  where $W_z$ is as in \lcite{\ComputingFT}.  The integrated form of
$\rho$ will also be denoted by $\rho$.  Since we know that $C^*(\B)$
is isomorphic to $\K$, a simple algebra, any representation of it will
be faithful and hence so is $\rho$.  We shall then identify $C^*(\B)$
with its image under $\rho$.

Let $\eta$ be in $C_c(\B)$, so that $\eta(z) = (f_z,z)$, for $z\in
S^1$, where $z\mapsto f_z$ is a continuous $C(S^1)$--valued function
on $S^1$.  We write $f_z(w)$ as $f(z,w)$.  In the latter form $f$
represents a continuous function on $S^1\x S^1$.  Observe that
  $$
  \rho(\eta)\xi\calcat z=
  \int_{S^1} f(w,z)\xi(w z) \d w =
  \int_{S^1} f(w z\inv,z)\xi(w) \d w,
  $$
  for $\xi\in C_c(S^1)\subseteq\H$.  Therefore $\rho(\eta)$ is the
(classical) integral operator with kernel $k(z,w) = f(w z\inv,z)$.
Since the transformation $(z,w)\mapsto (w z\inv,z)$ is a homeomorphism
of $S^1\x S^1$, we have that $C_c(\B)$ (or rather its image under
$\rho$) consists exactly of the integral operators with continuous
symbols.

For example, if $P$ is the rank one projection onto the space spanned
by $\phi\in \H$ with $\[\phi\]_2=1$, we have that
  $$
  P\xi\calcat z =
  \<\xi,\phi\>\phi \calcat z =
  \int_{S^1} \phi(z)\overline{\phi(w)}\xi(w) \d w,
  $$
  so that $P$ is also an integral operator.  Its kernel, given by
$k(z,w)=\phi(z)\overline{\phi(w)}$, is clearly continuous if and only
if $\phi$ is a continuous function.  Therefore $P$ belongs to
$C_c(\B)$ if and only if $\phi\in C(S^1)$.

Fix a measurable function $\delta\in L_\infty(S^1)$ such that
$|\delta(z)|=1$ for all $z$.  The multiplication operator $M_\delta$
is then a unitary operator on $\H$ which commutes with the bilateral
shift and hence the map
  $$
  \Delta: T\in\K \mapsto M_\delta T M_\delta\inv\in \K
  $$
  is an automorphism of $\K$ which commutes with $\a$.  It follows
that all aspects of the dynamical system $(\K,\Z,\a)$ are left
invariant under $\Delta$.  For example, if $T$ is $\a$-integrable then
so is $\Delta(T)$ and
  $$
  \F {\Delta(T)}z = \Delta(\F Tz)
  \for z\in S^1,
  $$
  and so on.

If $T$ is an integral operator with kernel $k$, observe that
  $$
  \Delta(T)\xi\calcat z =
  M_\delta T M_\delta^* \xi \calcat z =
  \int_{S^1} \delta(z) k(z,w) \overline{\delta(w) }\xi(w) \d w
  \for \xi\in C_c(S^1),
  $$
  and hence that $\Delta(T)$ is the integral operator with kernel
  $$
  k'(z,w) = \delta(z) k(z,w) \overline{\delta(w)}.
  $$
  So we see that $\Delta(C_c(\B))$ consists exactly of the integral
operators whose kernels have the above form for a continuous function
$k$.

If $\delta$ is sufficiently discontinuous we may then have
$\Delta(C_c(\B))\neq C_c(\B)$.  Nevertheless, given that $\Delta$
commutes with $\a$, one can use $\Delta(C_c(\B))$ in place of
$C_c(\B)$ in \lcite{\MoritaForFellBundles} in order to obtain a {\MR}
equivalence between $B_e$ and an ideal of $K\crossproduct_\a \Z$.  Let
  $$
  \P_1 = C_c(\B)^2 \cap C^*(\B)_+
  \and
  \P_2 = \Delta(C_c(\B))^2 \cap C^*(\B)_+.
  $$
  We claim that $\P_1\cup\P_2$ is not necessarily relatively
continuous.  In fact, let $P$ be the rank-one projection onto the
space spanned by a unit vector $\phi$ (for the 2-norm) belonging to
$C(S^1)$.  Then it is easy to see that
  $\Delta(P)$ is the projection onto the span of $\delta\phi$.  Let
$Q$ be another rank-one projection whose range is spanned by a unit
vector $\psi\in C(S^1)$, then as seen above
  \iItem{$\bullet$} $P\rc Q$ because $\overline\phi\psi \in C(S_1)$,
  \Item{$\bullet$} $\Delta(P)\rc \Delta(Q)$ because
$\overline{\delta\phi}\delta\psi = \overline\phi\psi \in C(S_1)$.
  \medskip\noindent
  However, it may happen that $P\not\rc \Delta(Q)$ because nothing
guarantees that $\overline{\phi}\delta\psi \in C(S_1)$.

One of the main lessons to be learned from this example is that a
maximal relatively continuous cone $\P$, as discussed in
\lcite{\PropertiesOfMaximalCones} and \lcite{\AbstractMR}, is not
unique.  In fact, by Zorn's Lemma we may take for each $i=1,2$, a
maximal such cone containing $\P_i$, say $\widetilde\P_i$, and we must
then have $\widetilde\P_1 \neq \widetilde\P_2$, or else $\P_1\cup\P_2$
is relatively continuous.

  This shows that a choice has to be made somewhere if question
\lcite{\DensityOfRCCone} is to be answered affirmatively.  Clearly the
same goes for Question \lcite{\MainQuestion}.

Let $\X_i = \overline{\zeta(\N_i)}$, for $i=1,2$, where $\N_i$ is
constructed from $\P_i$ as in \lcite{\AbstractMR}.  It is easy to see
that $\overline{\X_1^*\X_1}=\overline{\X_2^*\X_2}$ so that we fall
short of giving a counter-example for Question \lcite{\QuestionOnUFA}
even if we knew (which we don't) that the $\P_i$ were maximal.

  \endgroup

  \bye